\documentclass{elsart}
\usepackage{graphicx}
\usepackage{latexsym}
\usepackage{amssymb}
\usepackage{subfigure}
\usepackage{epsfig}

\newcommand{\PM}{Petviashvili }

\begin{document}
\begin{frontmatter}
\title{Petviashvili type methods for traveling wave computations: I. Analysis of convergence}
\author{J. \'Alvarez}
\address{Department of Applied Mathematics,
University of Valladolid, Paseo del Cauce 59, 47011,
Valladolid, Spain.}
\address{
IMUVA, Institute of Mathematics of University of Valladolid; Spain.
Email: joralv@eii.uva.es}
\author{A. Dur\'an \thanksref{au}}
\address{Department of Applied Mathematics, University of
Valladolid, Paseo de Bel\'en 15, 47011-Valladolid, Spain.}
\address{
IMUVA, Institute of Mathematics of University of Valladolid; Spain.
Email:
angel@mac.uva.es }

\thanks[au]{Corresponding author}

\begin{abstract}
In this paper a family of fixed point algorithms for the numerical resolution of some systems of nonlinear equations is designed and analyzed. The family introduced here generalizes the Petviashvili method and can be applied to the numerical generation of traveling waves in some nonlinear dispersive systems. Conditions for the local convergence are derived and numerical comparisons between different elements of the family are carried out.
\end{abstract}
\begin{keyword}
Petviashvili type methods, traveling wave generation, iterative
methods for nonlinear systems, orbital convergence

MSC2010: 65H10, 65M99, 35C99, 35C07, 76B25
\end{keyword}
\end{frontmatter}

\section{Introduction}
Considered here is the construction and study of fixed point algorithms for the numerical resolution of nonlinear systems of the form
\begin{eqnarray}
L{ u}=N({ u}),\quad u\in \mathbb{R}^{m}, \quad m>1,\label{m1}
\end{eqnarray}
where $L$ is a nonsingular $m\times m$ real matrix and
$N:\mathbb{R}^{m}\rightarrow \mathbb{R}^{m}$ is an homogeneous
function of the components of $u$ with degree $p, |p|>1$. These systems are typical in many applications, including the approximation to equilibria in mechanical systems and the numerical generation of traveling waves and ground states in nonlinear dispersive systems for water waves and nonlinear optics. {(In this last context, $m$ would represent the number of discretization points.)} More generally, (\ref{m1}) may appear when generating relative equilibria or coherent structures, \cite{champneyss}. {We denote by $u^{\ast}$ a solution of (\ref{m1}), that is
\begin{eqnarray}
Lu^{\ast}=N(u^{\ast}).\label{m1a}
\end{eqnarray}
}
The classical fixed point algorithm for (\ref{m1}) has the following formulation. If $u_{0}\neq 0$, the approximation to  $u^{*}$ in (\ref{m1a}) at the $(n+1)$-th iteration is given by the recurrence
\begin{eqnarray}
  Lu_{n+1}=N(u_{n}), \quad n=0,1,\ldots\label{m1b}
\end{eqnarray}

{ The method (\ref{m1b}) is not usually convergent for this kind of problems. Note that if
\begin{eqnarray}
S=L^{-1}N^{\prime}(u^{\ast}),\label{m2b}
\end{eqnarray}
 stands for the iteration matrix at $u^{\ast}$ (and where $N^{\prime}(u)$ denotes the Jacobian of $N$ at $u$), then the homogeneous character of $N$ implies that
$
N^{\prime}(u^{\ast})u^{\ast}=pN(u^{\ast});
$ therefore, using (\ref{m1a}),
\begin{eqnarray*}
S(u^{\ast})u^{\ast}=L^{-1}N^{\prime}(u^{\ast})u^{\ast}=pL^{-1}N(u^{\ast})=pu^{\ast}.
\end{eqnarray*}
Thus, $u^{\ast}$ is an eigenvector of $S$ associated to an eigenvalue $\lambda=p$ with $|p|>1$.
}

The methods presented here generalize the so-called Petviashvili method. From a starting iteration $u_{0}\neq 0$, the \PM method generates the recurrence
\begin{eqnarray}
m(u_{n})&=&\frac{\langle Lu_{n},u_{n}\rangle}{\langle N(u_{n}),u_{n}\rangle},\label{m2}\\
  Lu_{n+1}&=&m(u_{n})^{\gamma}N(u_{n}), \quad n=0,1,\ldots\label{m3}
\end{eqnarray}
where {here and in the rest of the paper} $\langle\cdot,\cdot\rangle$ stands for the Euclidean inner
product and $\gamma$ is a free real parameter. The term (\ref{m2})
is called stabilizing factor and, in the case of convergence, must
tend to one.
The origin of the method is in \cite{petviashvili}, focused on the
search for lump solitary waves of the Kadomtsev-Petviashvili I
(KPI) equation
\begin{eqnarray}
\left(u_{t}+2uu_{x}+u_{xxx}\right)_{x}=u_{yy},\quad t>0, x,y\in \mathbb{R}\label{crete1}
\end{eqnarray}
of the form $u(x,y,t)=c\varphi(X,Y)=c\varphi(\sqrt{c}(x-ct),cy), c>0$.
The profile $\varphi$ must satisfy
\begin{eqnarray*}
\partial_{XX}\left(-\varphi+\partial_{XX}\varphi\right)-\partial_{YY}\varphi=
-\partial_{XX}\varphi^{2},
\end{eqnarray*}
which, in terms of the
$2$-D Fourier Transform,
\begin{eqnarray*}
\widehat{\varphi}({k_{x}, k_{y}})=\int_{-\infty}^{\infty}\int_{-\infty}^{\infty} \varphi(x,y)e^{-ik_{x}x}e^{-ik_{y}y}dxdy,
\end{eqnarray*}
is converted into an algebraic system
\begin{eqnarray}
&&\widehat{\varphi}({k_{x}, k_{y}})=G({k_{x}, k_{y}})A({k_{x}, k_{y}}),\nonumber\\
&&G({k_{x}, k_{y}})=\frac{k_{x}^{2}}{k_{x}^{4}+k_{x}^{2}+k_{y}^{2}},\quad
A({k_{x}, k_{y}})=\widehat{\varphi^{2}}({k_{x}, k_{y}})\label{crete3}
\end{eqnarray}
 The divergence of the classical fixed point algorithm, applied to (\ref{crete3})
 forces to consider a new iteration system, of the form
\begin{eqnarray*}
\widehat{\varphi}({k_{x}, k_{y}})=m(\varphi)^{\gamma}G({k_{x}, k_{y}})A({k_{x}, k_{y}}),
\end{eqnarray*}
where the stabilizing factor $m(\varphi)$ is defined as
\begin{eqnarray*}
m(\varphi)=\frac{s_{1}}{s_{2}},\quad
s_{1}=\int\int|\widehat{\varphi}|^{2}dk_{x}dk_{y},\,
s_{2}=\int\int\frac{k_{x}^{2}}{k_{x}^{4}+k_{x}^{2}+k_{y}^{2}}\widehat{\varphi^{2}}
\overline{\widehat{\varphi}}dk_{x}dk_{y},
\end{eqnarray*}
{($\overline{\widehat{\varphi}}$ denotes the complex conjugate of ${\widehat{\varphi}}$)} where $\gamma$ is a free real parameter. In the case of
(\ref{crete1}), $\gamma$ is taken approximately $2$, \cite{petviashvili}. On the other hand, for the
exact profile $\varphi$, $m(\varphi)=1$.

The \PM method has become popular as a technique to generate special solutions in partial differential equations of interest in water waves and nonlinear optics. It takes part of a large family of methods designed to this goal, which includes variants of the Newton's method, \cite{yang}, {modified conjugate gradient methods applied to nonlinear problems, \cite{lakoba}}, squared operator methods, \cite{yangl2},
imaginary-time evolution methods, \cite{yangl1} or
different variational procedures, \cite{garciap,baod,caliario} .
Some literature about (\ref{m2}), (\ref{m3}), from the original paper, \cite{petviashvili}, is now briefly
reviewed. Pelinovsky and Stepanyants, \cite{pelinovskys}, analyze
the continuous version of the method to approximate solitary wave
profiles of the nonlinear dispersive models
\begin{eqnarray*}
    u_{t}-\mathcal{L}u_{x}+pu^{p-1}u_{x}&=&0,\quad p>1, \quad t>0, \quad x\in\mathbb{R},
\end{eqnarray*}
where $\mathcal{L}$ is a pseudodifferential operator with positive
Fourier symbol. On the other hand, Lakoba and Yang, \cite{lakobay,lakobay2}, introduce a generalized version of the procedure, for more general systems of the form
\begin{eqnarray*}
-Mu+F(x,u)=0,\quad u\rightarrow 0, \quad |x|\rightarrow\infty
    \end{eqnarray*}
    where $M$ is positive definite, self-adjoint operator and $F$ is nonlinear (see also \cite{yang2}).
Finally,
Ablowitz and Musslimani,
\cite{ablowitzm}, (see also \cite{ablowitzm2}) propose an alternative of the algorithm, the spectral renormalization method, with
application to generate numerically ground state profiles for systems of NLS
type
\begin{eqnarray*}
    iU_{z}+\Delta U-V(x)U+f(|U|^{2})U=0.
    \end{eqnarray*}

Some new results contained in this paper are described below.


\begin{itemize}
\item
Based on the philosophy the \PM method was devised with, new fixed point methods are derived. They can be considered as a Petviashvili type family of methods.
\item From the view point of the classical algorithm, the corresponding iteration functions are designed to filter the harmful directions of the errors leading to divergence, in such a way that convergence results are obtained under the same hypotheses as those of the \PM method. Here it is worth mentioning two types of convergence. The first one has the classical sense, with the requirement (among others) of isolated fixed points. However, in traveling wave generation, it is very typical that the system of equations admits a symmetry group {(usually related to translational or rotational invariance of the system)}. In this case, fixed points cannot be isolated and the convergence must be understood in the orbital sense, that is, for the orbits of fixed points. Convergence results for both cases will be given. This study complements some previous results of convergence presented in the literature for the \PM method, \cite{pelinovskys,lakobay,lakobay2}.
\end{itemize}
The structure of the paper is as follows: in Section \ref{sec2} and starting from the \PM method (\ref{m2}), (\ref{m3}), the new family of fixed point algorithms are constructed and analyzed. A comparison of efficiency of some of them is also carried out. Section \ref{sec3} will treat the
derivation of general conditions for the local convergence of the
methods. The first part of the study assumes the existence of a neighborhood where the fixed point $u^{*}$ is unique. The spectral analysis of the iteration matrix (\ref{m2b}) of the classical
fixed point algorithm (or, equivalently, the pencil $A(\lambda)=\lambda L-N^{\prime}(u^{*})$) is used.
The local convergence of the methods can be
achieved even when (\ref{m2b}) admits eigenvalues with modulus greater
than or equals one. The results are illustrated with several numerical examples, concerning the generation of localized ground state solutions of nonlinear Schr\"{o}dinger type equations with potentials. On the other hand, some other applications of the methods suggest to
analyze a case where the hypothesis of local uniqueness of the fixed point does not hold, in the sense that the system (\ref{m1}) admits a group of symmetries, generating orbits of solutions. From the point of view of the analysis, the existence of a symmetry group in (\ref{m1}) is associated to the formation of the eigenvalue one in the pencil, \cite{champneyss}, and leads, in a natural way, to the concept of orbital convergence.
Section \ref{sec3} is finished off with the corresponding results of convergence for this case and they will be illustrated by the generation of soliton solutions of the nonlinear Schr\"{o}dinger equation.

The present paper is a first part of a study of the methods carried out by the same authors. It will be followed by a second part, in which some particular, relevant cases of systems (\ref{m1}) are emphasized and where the effect of the introduction of acceleration techniques is studied.

%

\section{Derivation and convergence analysis of the algorithms}
\label{sec2}
\subsection{Derivation}
The following fixed point methods for the iterative resolution of (\ref{m1}) are introduced. If ${ u}_{0}\neq 0$, the iterations ${u}_{n}, n=1,2,\ldots$ are generated by a formula of the form
\begin{eqnarray}
 Lu_{n+1}=s(u_{n})N(u_{n}), \quad n=0,1,\ldots\label{m4}
\end{eqnarray}
where $s:\mathbb{R}^{m}\rightarrow \mathbb{R}$ is a $C^{1}$ function satisfying the following properties:
\begin{itemize}
\item[(P1)] A set of fixed points of the iteration operator
\begin{eqnarray}
F(u)=s(u)L^{-1}N(u),\label{iterop}
\end{eqnarray}
coincides with a set of fixed points of (\ref{m1}). This means that: (a) if $u^{*}$ is a solution of (\ref{m1}) then $s(u^{*})=1$; (b) inversely, if the sequence $\{u_{n}\}_{n=0}^{\infty}$, generated by (\ref{m4}), converges to some $y$, then $s(y)=1$ (and, consequently, $y$ is a solution of (\ref{m1})).
\item[(P2)] $s$ is homogeneous with degree $q$ such that $|p+q|<1$.
\end{itemize}
Note that, in particular, the choice
\begin{eqnarray}
s(u)=\left(\frac{\langle Lu,u\rangle}{\langle N(u),u\rangle}\right)^{\gamma},\quad q=\gamma(1-p),\label{m5}
\end{eqnarray}
leads to the \PM method (\ref{m2}), (\ref{m3}). Then, (\ref{m4}) can be considered as a generalization and justifies that $s$ will be also called a stabilizing factor. Several examples are the following:
\begin{itemize}
\item The term (\ref{m5}) can be generalized by considering any $C^{1}$ homogeneous function $f:\mathbb{R}^{m}\rightarrow \mathbb{R}$ with degree greater than or equals one and taking
 \begin{eqnarray}
s_{f}(u)=\left(\frac{\langle Lu,f(u)\rangle}{\langle N(u),f(u)\rangle}\right)^{\gamma},\quad q=\gamma(1-p), \quad |p+q|<1.\label{m6a}
\end{eqnarray}
\item Another alternative is the use of norms, with
\begin{eqnarray}
s_{r}(u)=\left(\frac{||Lu||_{r}}{||N(u)||_{r}}\right)^{\gamma},\quad q=\gamma(1-p), \quad |p+q|<1,\label{m6b}
\end{eqnarray}
where if $u=(u_{1},\ldots,u_{m})^{T}$ then $||u||_{r}=\left(|u_{1}|^{r}+\ldots+|u_{m}|^{r}\right)^{1/r}, 1\leq r\leq +\infty$, with $r=+\infty$ standing for the usual maximum norm. The case $r=1$ was considered in \cite{ablowitzm2}.
\end{itemize}

\subsection{First numerical experiments}
\label{sec22}
Displayed here are some numerical experiments concerning the performance of the methods, according to the choice of the stabilizing factor. As an example
we consider the problem of generating lump solitary
waves in the 2D
Benjamin equation
\begin{eqnarray}
\label{ben2da}
\left(\eta_{t}+\alpha (\eta^{2})_{x}-\beta \mathcal{H}(\eta_{xx})+\delta \eta_{xxx} \right)_{x}-\eta_{zz}=0,
\end{eqnarray}
where $\alpha, \beta, \delta\geq 0$ and $\mathcal{H}$ stands for the Hilbert transform with respect to $x$:
\begin{eqnarray*}
\mathcal{H}f(x)=\frac{1}{\pi}P.V. \int_{-\infty}^{\infty}
\frac{f(y)}{x-y}dy.
\end{eqnarray*}
Equation (\ref{ben2da}) is analyzed in \cite{kim,kima1,kima2}. It
appears as an extension of the one-dimensional equation derived by
Benjamin, \cite{ben0,ben1,ben2}, and modeling the propagation of
waves at the interface of two ideal fluids, with a bounded upper
layer and the heavier one with infinite depth, and under the
presence of interfacial tension. The two-dimensional version
incorporates weak transverse variations. The form (\ref{ben2da})
contains particular cases, such as the Kamdotsev-Petviashvili
(KP-I) equation, \cite{kadomtsevp,manakovzbim} ($\beta=0,
\delta>0$) and the two-dimensional Benjamin-Davis-Ono (BDO)
equation, corresponding to $\delta=0, \beta>0$, \cite{ablowitzs}. A
normalized form of (\ref{ben2da})
\begin{eqnarray}
\label{ben2db}
\left(\eta_{t}+ (\eta^{2})_{x}-2{\Gamma} \mathcal{H}(\eta_{xx})+\eta_{xxx} \right)_{x}-\eta_{zz}=0,
\end{eqnarray}
is derived in \cite{kima2} and will be adopted here. The parameter
$\Gamma\geq 0$ is related to the interfacial tension and the
densities of the fluids. Finally, for localized solutions, the
constraint
\begin{eqnarray}
\label{constraint}
\int_{-\infty}^{\infty} \eta(x,z,t)dx=0,
\end{eqnarray}
(zero total mass condition) is assumed, as in the KP and BDO
equations \cite{katsisa}.
\begin{figure}[htbp]
\centering \subfigure{\subfigure{
\includegraphics[width=5.5cm]{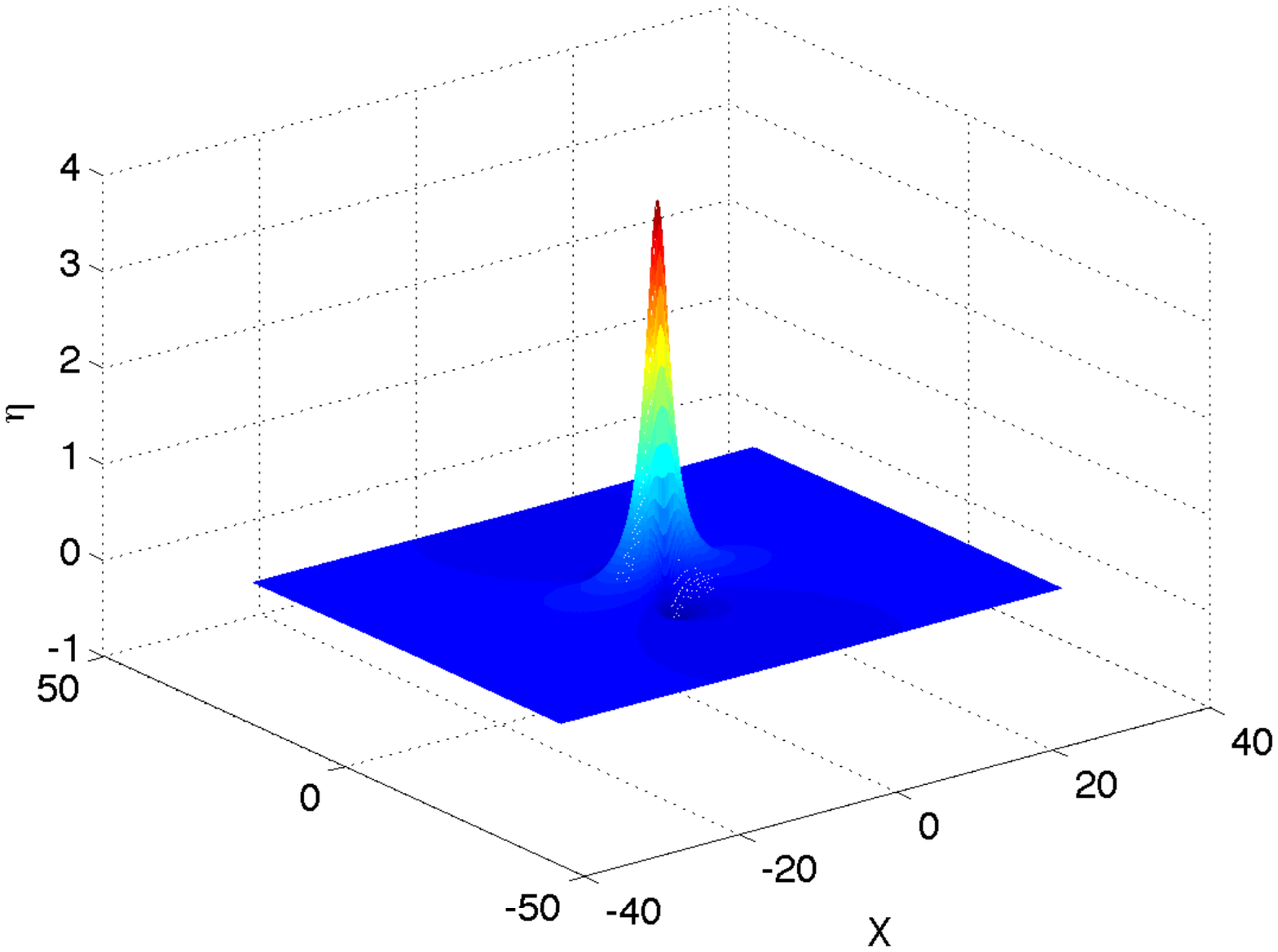}}
\subfigure{
\includegraphics[width=5.5cm]{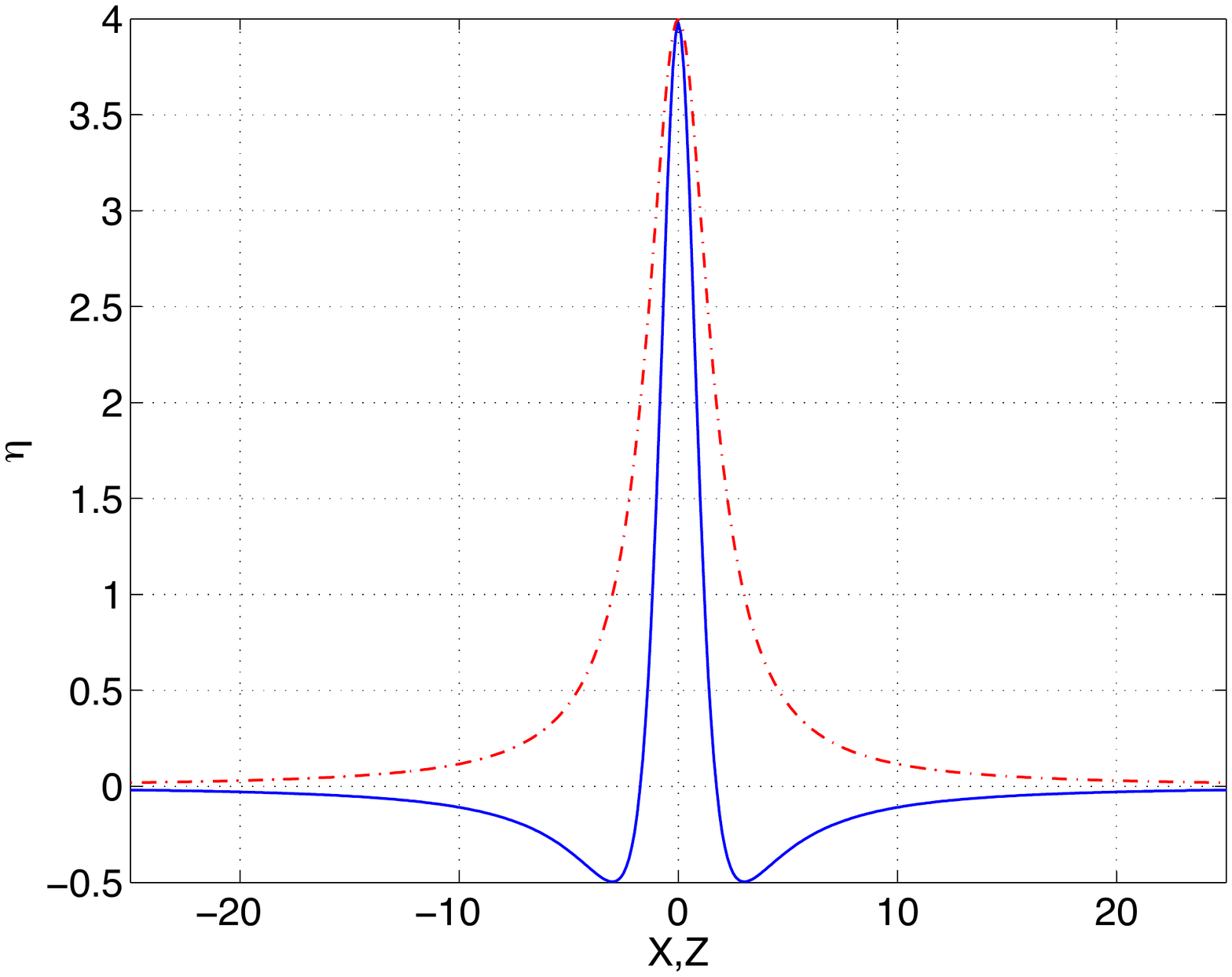}}
}
\end{figure}
\begin{figure}[htbp]
\centering \subfigure{\subfigure{
\includegraphics[width=5.5cm]{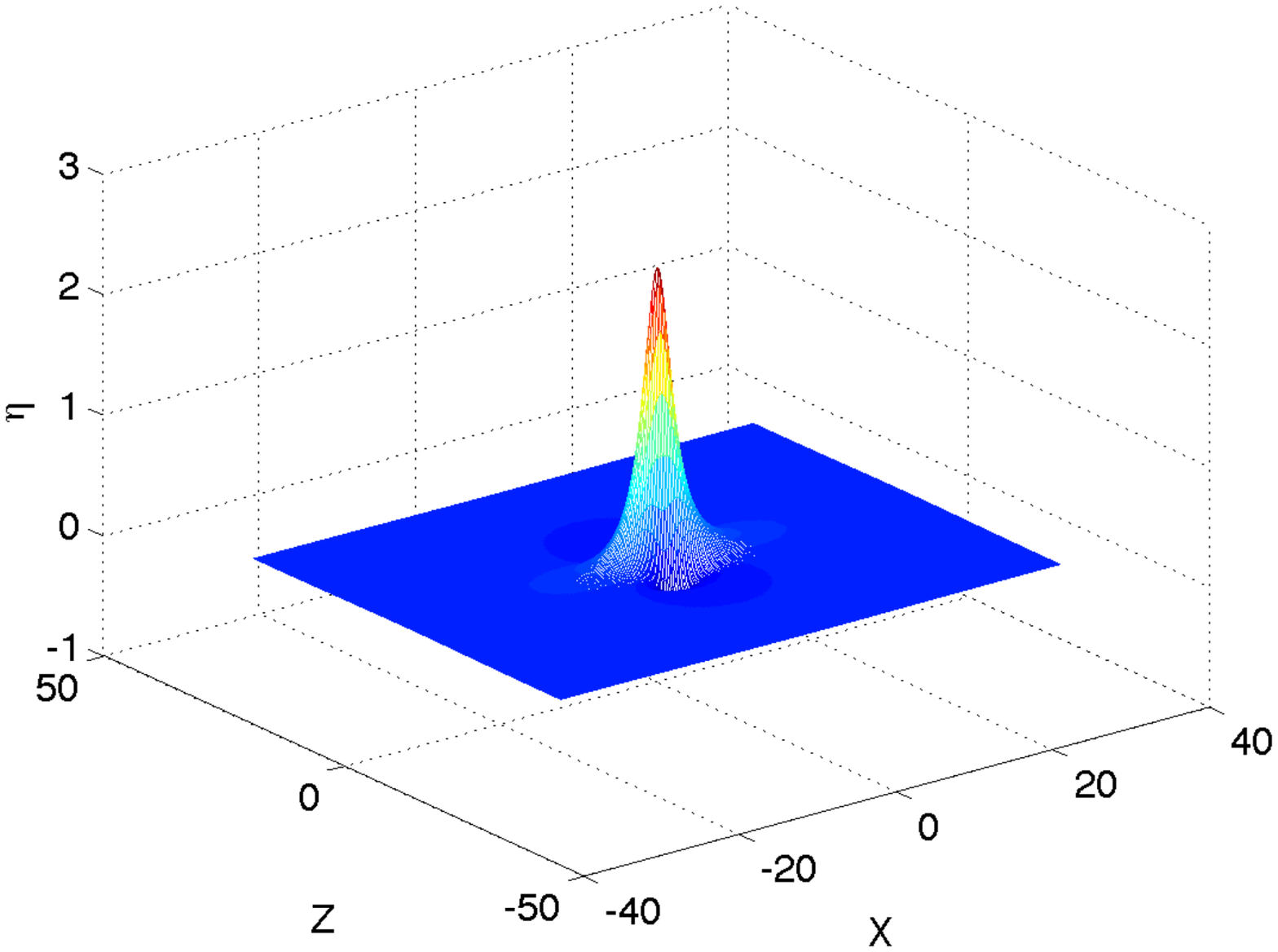}}
\subfigure{
\includegraphics[width=5.5cm]{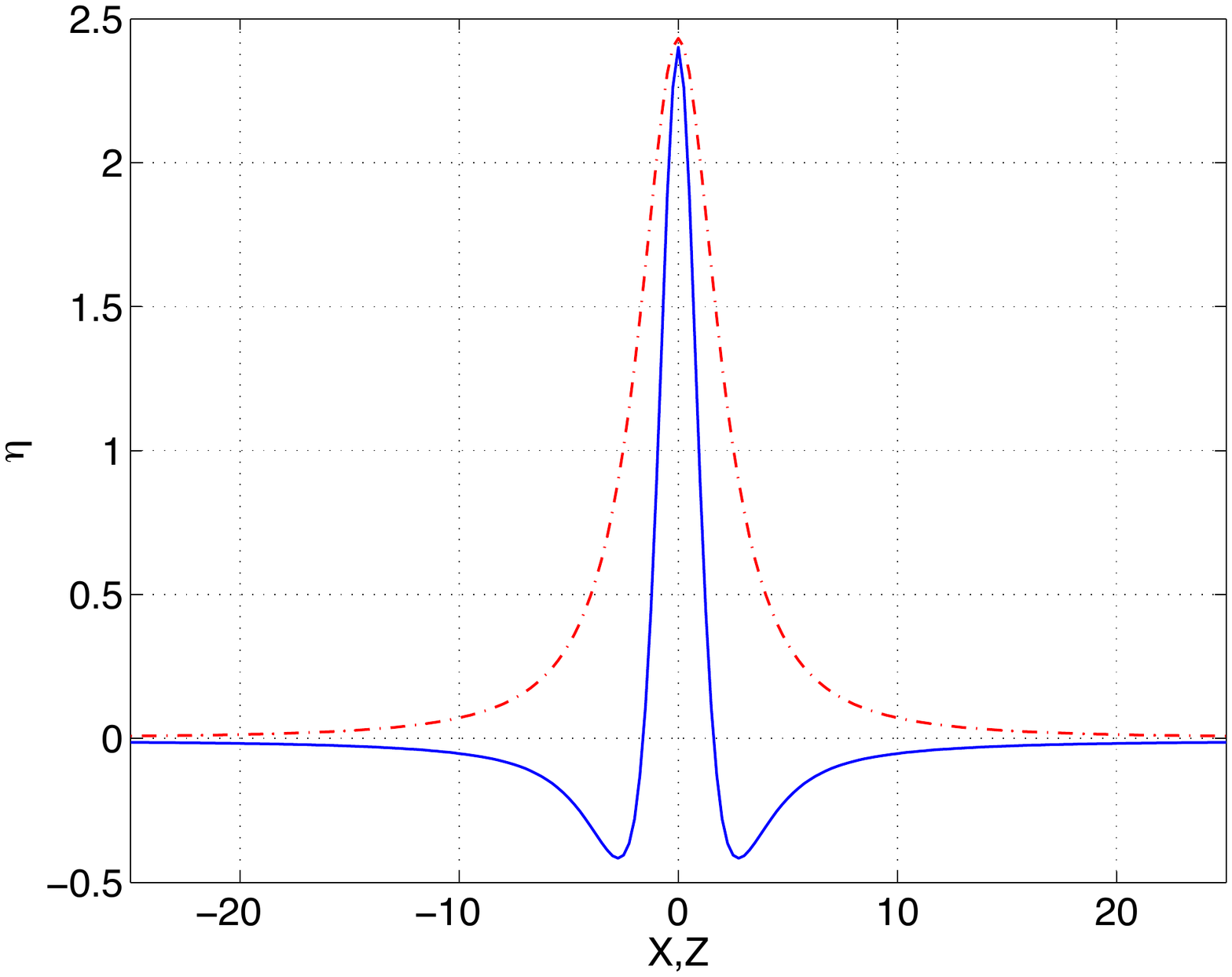}} }
\end{figure}
\begin{figure}[htbp]
\centering \subfigure{\subfigure{
\includegraphics[width=5.5cm]{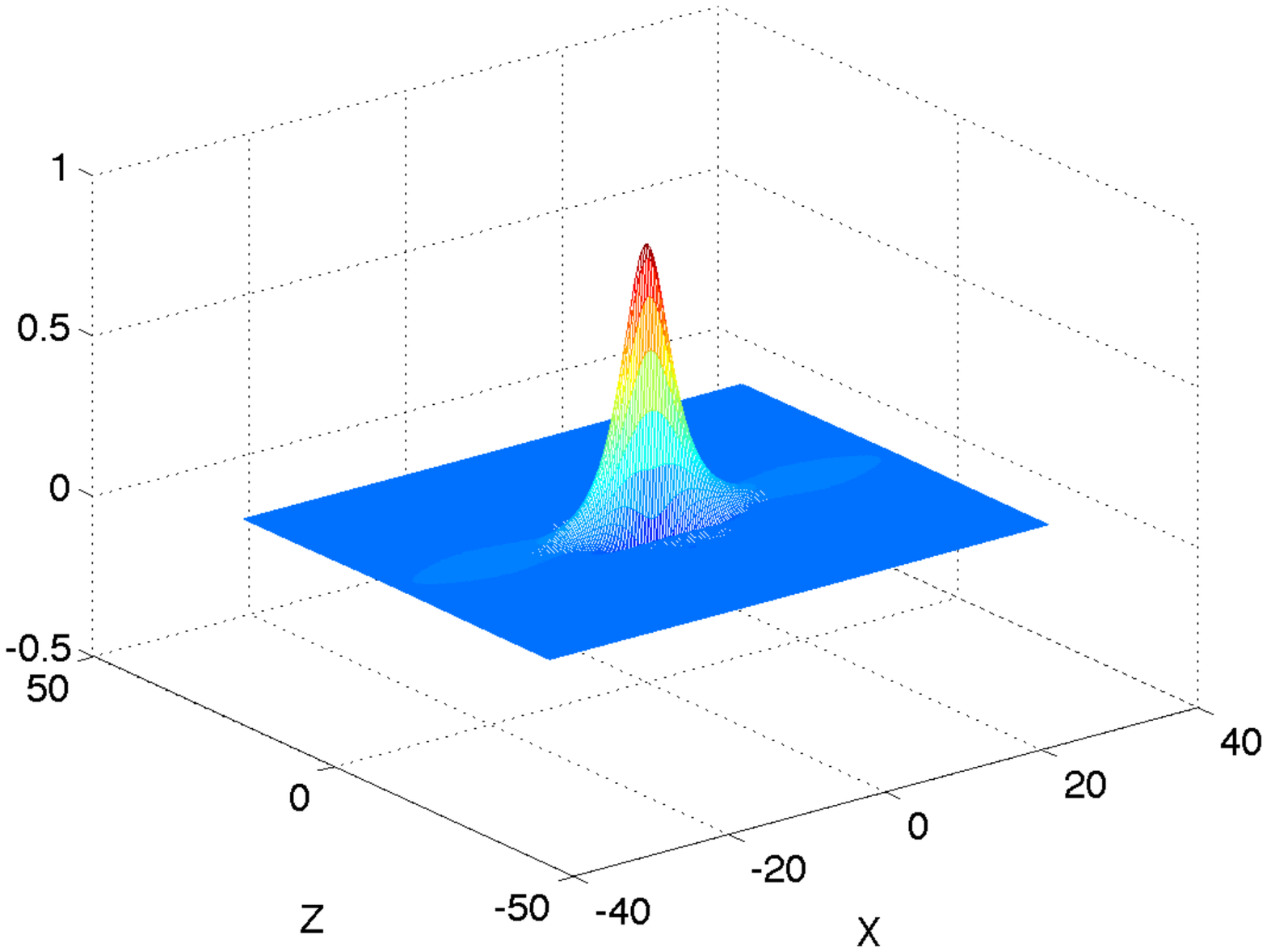}}
\subfigure{
\includegraphics[width=5.5cm]{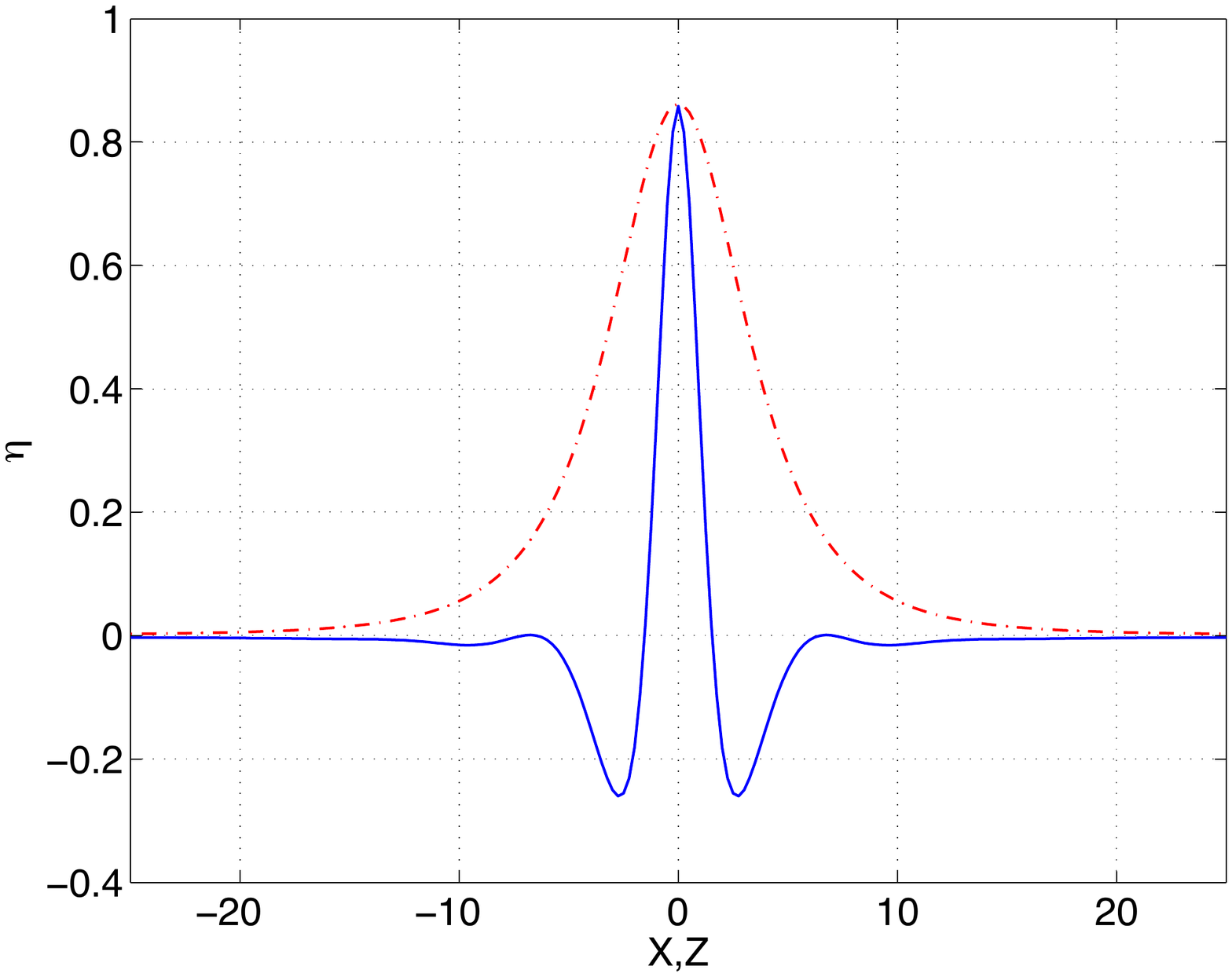}} }
\end{figure}
\begin{figure}[htbp]
\centering \subfigure{\subfigure{
\includegraphics[width=5.5cm]{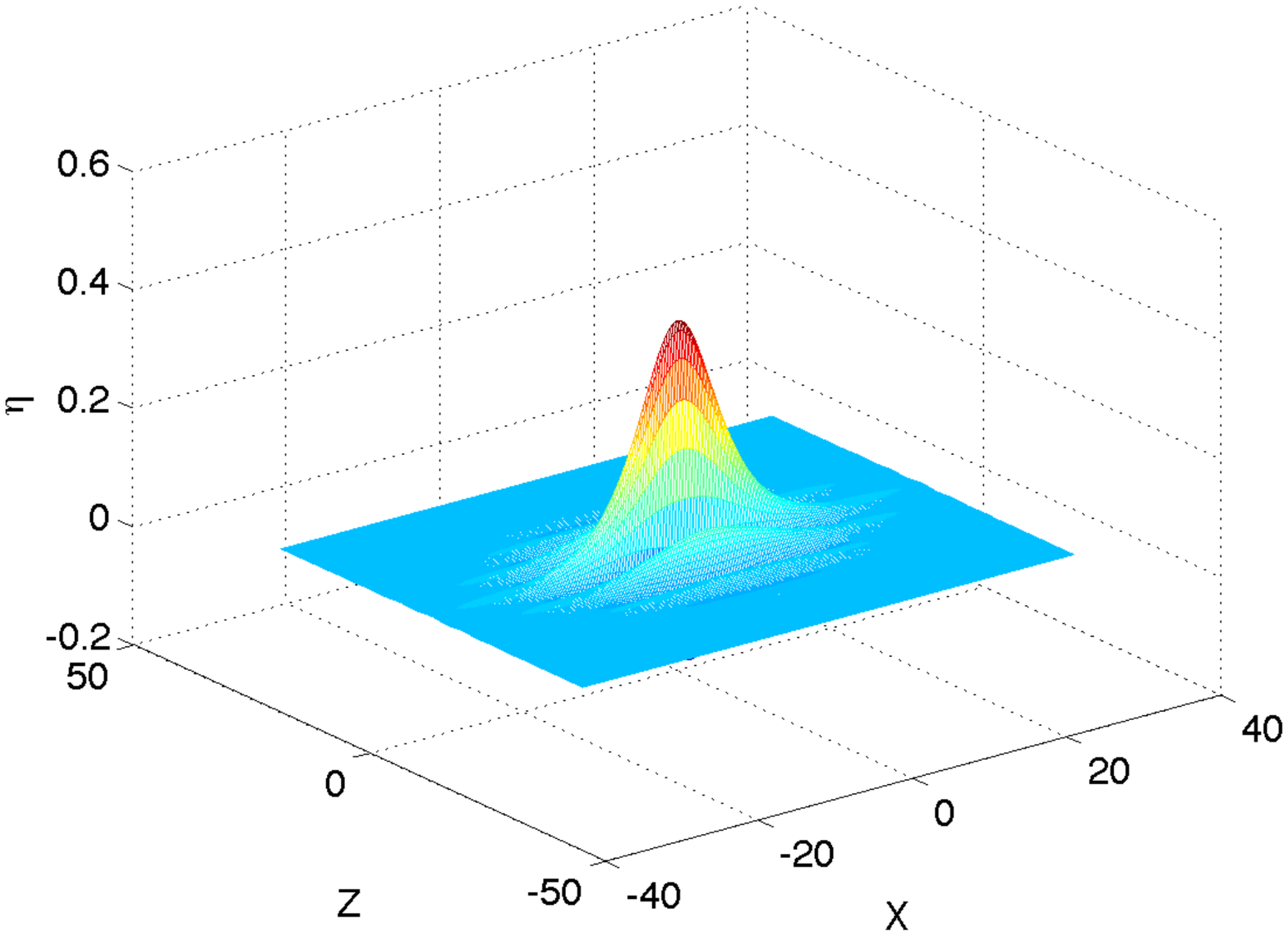}}
\subfigure{
\includegraphics[width=5.5cm]{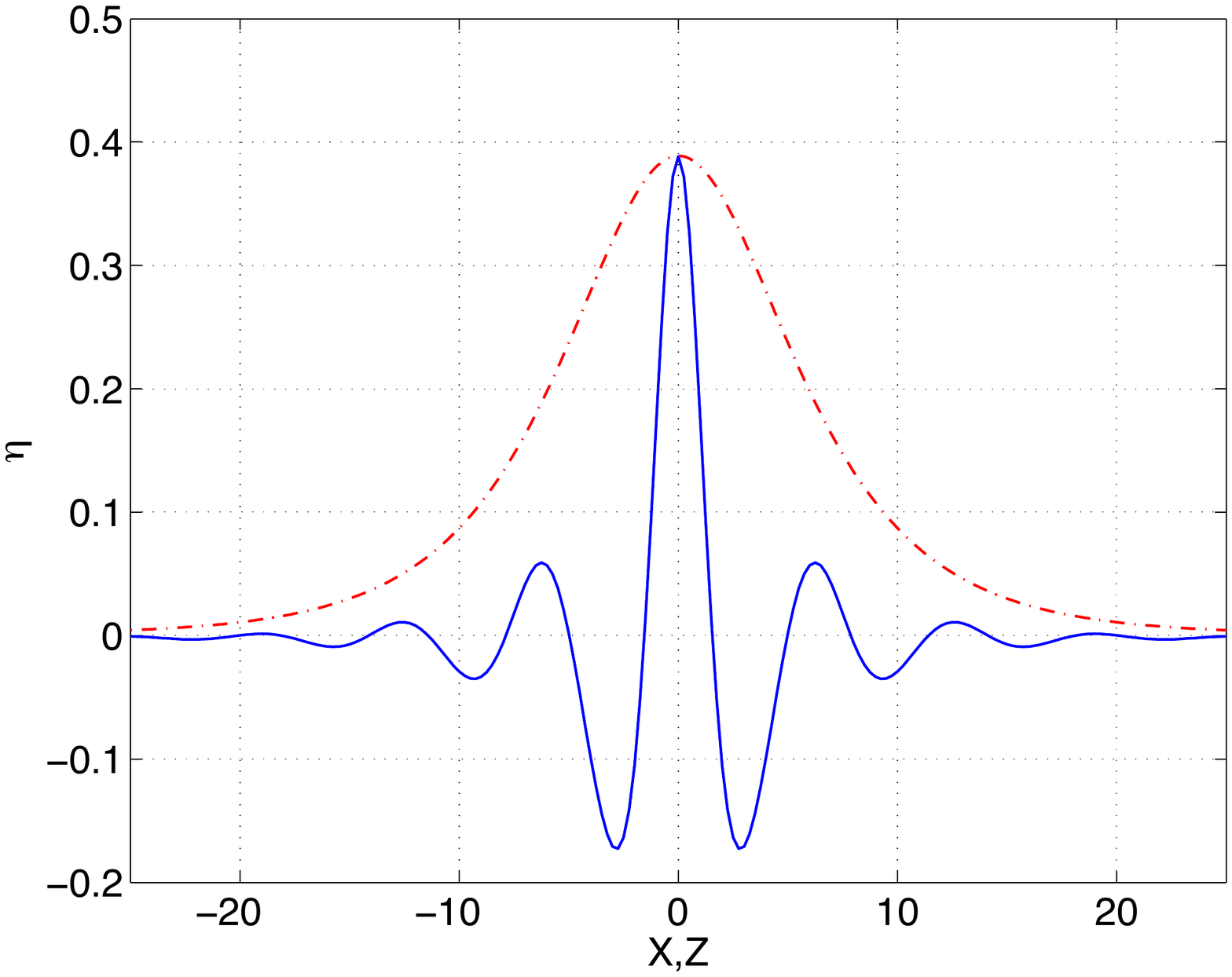}} }
\caption{Solitary wave generation of (\ref{ben2da}) with \PM
method. The approximated profiles correspond to $\Gamma=0.5,0.9,0.95,0.99$ (left).
On the right, the corresponding $X$ and $Z$ cross sections {are shown} (solid
and dashed-dotted lines, respectively) .} \label{fexample41}
\end{figure}

\begin{figure}[htbp]
\centering \subfigure[]{
\includegraphics[width=7cm]{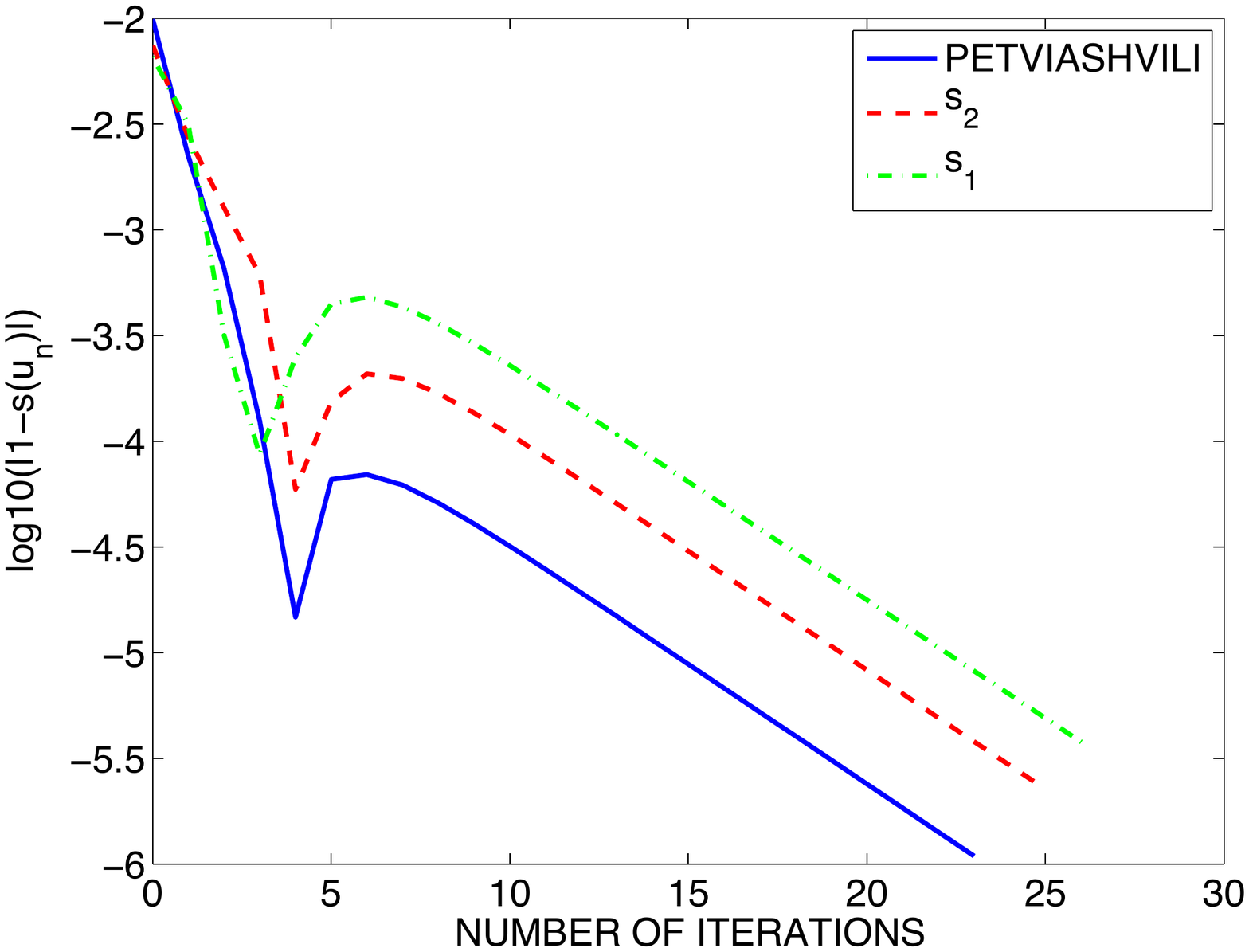} }
\subfigure[]{
\includegraphics[width=7cm]{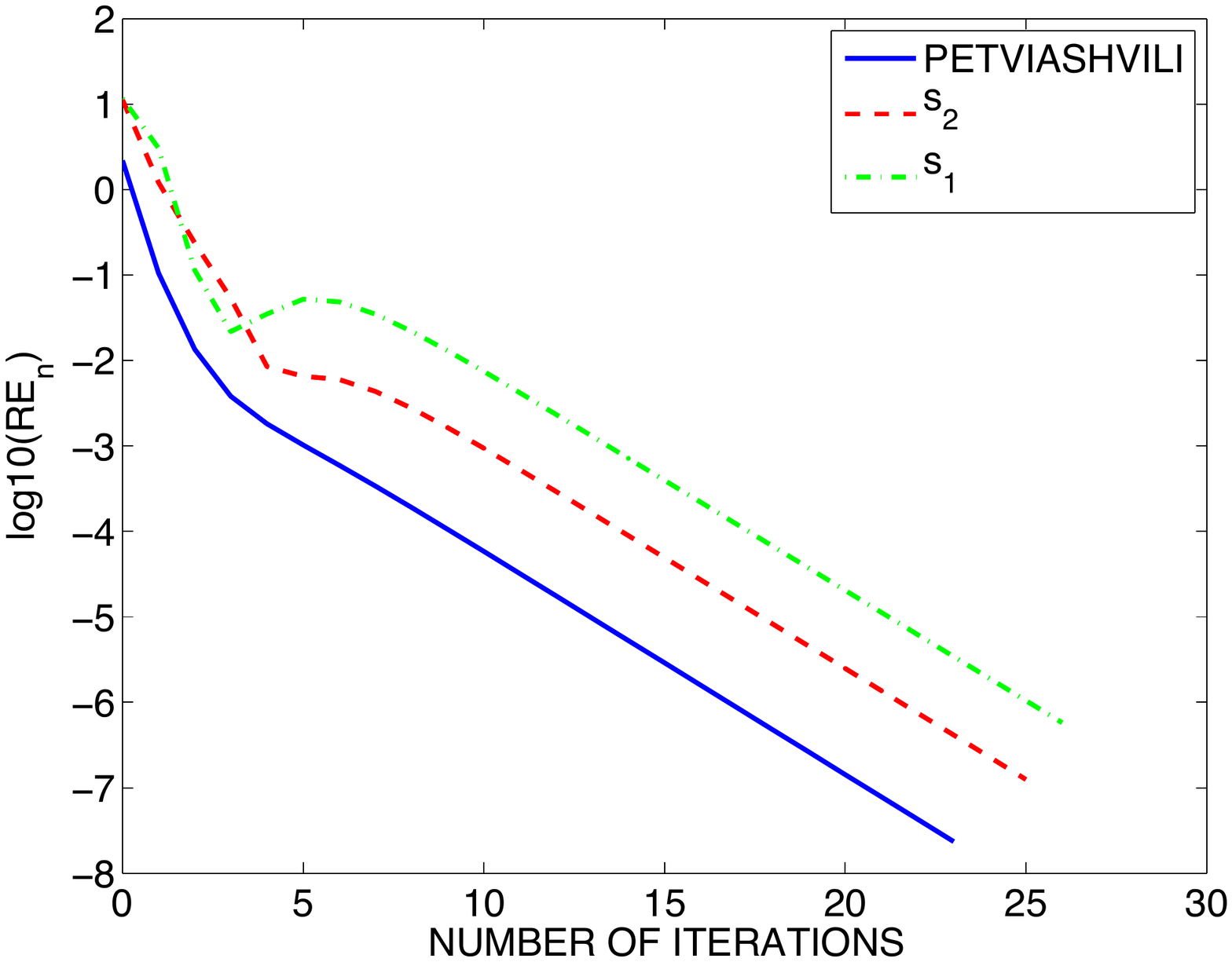} }
\caption{Convergence results of the \PM type methods (\ref{m4}) for (\ref{lumpsw}):
(a) Discrepancy between the stabilizing factor vs number of
iterations (semilog scale). (b) Logarithm of the residual errors vs number of
iterations. Solid line: Petviashvili method (\ref{m6a}) with {$f(x)=x$}; dashed line: (\ref{m6b}) with $r=2$; dashed-dotted line: (\ref{m6b}) with $r=1$.} \label{fexample42}
\end{figure}

The search for lump solitary wave solutions of (\ref{ben2db})
\begin{eqnarray*}
\eta(x,z,t)=\eta(X,Z),\quad X=x-c_{s}t, \quad Z=z,
\end{eqnarray*}
leads to the equation
\begin{eqnarray}
\label{lumpsw} \left(-c_{s}\eta+\eta^{2}-2\gamma
\mathcal{H}(\eta_{X})+\eta_{XX}\right)_{XX}-\eta_{ZZ}=0,
\end{eqnarray}
for the profile $\eta$. In \cite{kima2} it is shown that
(\ref{lumpsw}) admits lumps of wavepacket type, as well as lumps
of KP-I type {(see Figures \ref{fexample41}(d) and \ref{fexample41}(a) respectively)}. The bifurcation point corresponds to $\Gamma=1$. As
an alternative to the numerical procedure performed in that paper,
some lumps will be here generated numerically, by using numerical
continuation in $\Gamma$ from the known lump solitary wave
solution of the KP-I equation (corresponding to $\Gamma=0$),
\cite{manakovzbim}. Taking two-dimensional Fourier transform in
(\ref{lumpsw}) we have
\begin{eqnarray}
\left({k_{x}}^{2}\left(c_{s}+2\Gamma |k_{x}|+k_{x}^{2}\right)+{k_{y}}^{2}\right)\widehat{\eta}(k_{x},k_{y})={k_{x}^{2}}\widehat{\left(\eta^{2}\right)}(k_{x},k_{y}).\label{lumpsw2}
\end{eqnarray}
The continuation algorithm until the computation of the profile at certain value $\Gamma^{*}$ consists of defining an homotopic path
\begin{eqnarray*}
\Gamma_{0}=0<\Gamma_{1}<\cdots <\Gamma_{M}=\Gamma^{*}<1,
\end{eqnarray*}
and solving numerically (\ref{lumpsw2}) with the method (\ref{m4}) at
each $\Gamma_{j}$ and initial iteration given by the last
computed iterate at the previous stage $\Gamma_{j-1}$.
The procedure starts with the exact KP-I lump at
$\Gamma=\Gamma_{0}=0$. {The numerical resolution  of (\ref{lumpsw2})
with (\ref{m4}) is now described. Note that (\ref{lumpsw2}) for $k_{x}=k_{y}=0$ leaves the $(0,0)$-Fourier component free and the value $\widehat{\eta}(0,0)=0$
is set by the zero total mass condition (\ref{constraint}).  System (\ref{lumpsw2}) is discretized by using a Fourier collocation scheme and, with the notation of (\ref{m1}), the corresponding discrete system leads to a singular matrix $L$. In order to solve this, the $(0,0)$-Fourier component of the approximation is set to zero, then the resulting system for the rest of the Fourier components is not singular and is therefore iteratively solved with (\ref{m4}), \cite{pelinovskys}.}
%
%
The numerical results are shown in Figure
\ref{fexample41}. {(They correspond to $c_{s}=1$.)} On the left, different lump profiles, associated
to several values of $\Gamma$ are displayed. On the right, the
corresponding $X-$ and $Z-$ cross sections are represented. The
convergence of the procedure is illustrated by Figure
\ref{fexample42}. For the case $\Gamma=0.99$, it shows the limiting
behaviour of the stabilizing factor (left) and the residual errors
(right)
\begin{eqnarray}
RE_{n}=||Lx_{n}-N(x_{n})||,\label{res}
\end{eqnarray}
{(where here and in the rest of the paper $||\cdot ||$ stands for the Euclidean norm)}
both as functions of the number of iterations as for several choices of the stabilizing factor, from the two families (\ref{m6a}), (\ref{m6b}).
As a representative of (\ref{m6a}), the original Petviashvili method {($f(x)=x$)} has been compared with two methods with stabilizing factors of the form (\ref{m6b}), corresponding to $r=1,2$. Figure \ref{fexample42}(a) shows that the convergence of the stabilizing factor to one is more efficient with the Petviashvili method: it provides an error with less iterations and, for a fixed number of iterations, it gives a smaller discrepancy. This is also the conclusion when analyzing Figure \ref{fexample42}(b), concerning the behaviour of the residual error with the number of iterations. However, this better performance of the \PM method in this example is not big enough to be conclusive and to rule the rest of the methods out in a general situation. {We have the impression} that in terms of the computational effort, the methods are more or less equivalent, with a slight superiority of (\ref{m2}), (\ref{m3}). For that reason, this will be considered as a representative of (\ref{m4}) for the rest of the experiments in this paper (see Section \ref{sec3}).

{Finally, although the main goal of the example is illustrating a comparison between some methods of the family (\ref{m4}), it is worth mentioning that (\ref{ben2db}) is translationally invariant. Thus, the convergence must be understood in the sense analyzed in Section \ref{sec33} (orbital convergence).}
\section{Analysis of convergence}
\label{sec3}
As mentioned in the Introduction, the convergence of the methods (\ref{m4}) can be divided in two cases, depending on the character of $u^{*}$ as fixed point of (\ref{m1}).
{In what follows, the Jacobian of the iteration operator (\ref{iterop}) at a fixed point $u^{*}$ satisfying (\ref{m1a}),
\begin{eqnarray}
F^{\prime}(u^{*})=S+u^{*}\left(\nabla s(u^{*})\right),\label{m6c}
\end{eqnarray}
will be used. {(In (\ref{m6c}), the gradient $\nabla s(u^{*})$ is taken as a row vector.)}}
\subsection{Convergence (classical sense)}
We first define the pencil {$A(\lambda)=\lambda L-N^{\prime}(u^{*})$, where $u^{*}$ satisfies (\ref{m1a}). Note that, since $L$ is nonsingular, the zeros of $A(\lambda)$ coincide with the spectrum of the iteration matrix (\ref{m2b}), see \cite{Demmel}, Section 4.5 and \cite{GolubV}, Section 7.7. (In particular, $\lambda=p$ is a zero of $A(\lambda)$ with $A(p)u^{*}=0$.) We also remind that an eigenvalue $\lambda$ of a matrix is semisimple if the corresponding geometric and algebraic multiplicities are the same; that is, if the dimension of the associated eigenspace coincides with the order of $\lambda$ as zero of the characteristic polynomial.}
\begin{thm}
\label{theorem1}
Assume that
\begin{itemize}
\item[(H1)] There exists $R>0$ such that $u^{*}$ in (\ref{m1a}) is the unique fixed point of (\ref{m1})
in $B(u^{*},R)=\{u\in \mathbb{R}^{m} / ||u-u^{*}||<R\}$.
\end{itemize}
Take $u_{0}\neq 0$ and assume the following hypotheses on the zeros of $A(\lambda)$:
\begin{itemize}
\item[(i)] $\lambda=p$ is simple.
\item[(ii)] The rest of $\lambda$ satisfies $|\lambda|\leq 1$.
\item[(iii)] If $|\lambda|=1\Rightarrow\left\{\begin{matrix}\lambda \quad \mbox {is semisimple}\\
u_{0} \quad \mbox {does not have component in }
Ker A(\lambda)\end{matrix}\right.$
\end{itemize}
Then the method (\ref{m4}), with $s$ satisfying (P1) and (P2), is locally convergent, that is, there is a neighborhood $W$ of $u^{*}$ such that if $u_{0}\in W, u_{0}\neq 0$, the sequence $\{u_{n}\}_{n=0}^{\infty}$ generated by (\ref{m4}), converges to $u^{*}$. The optimal rate of convergence is obtained with $q=-p$.
\end{thm}

{
{\em Proof}.
The errors $e_{n}=u_{n}-u^{*}, n=0,1,\ldots,$ satisfy
\begin{eqnarray}
e_{n+1} =F^{\prime}(u^{*})e_{n}
+O(||e_{n}||^{2}),\quad n=0,1,\ldots,\label{m7}
\end{eqnarray}
where $F^{\prime}(u^{*})$ is given by (\ref{m6c}). According to the hypotheses (i)-(iii), $e_{n}$ can be decomposed
\begin{eqnarray}
e_{n}=\alpha_{n}u^{*}+z_{n}, \quad \alpha_{n}\in \mathbb{R}, \quad z_{n}\in
V, \quad S(V)\subset V,\label{m8}
\end{eqnarray}
where $V$ is a $S$ invariant supplementary subspace of $span(u^{*})$. ($V$ is the sum of the $S$-invariant subspaces associated to the eigenvalues of $S$ different from $\lambda=p$.) Substituting (\ref{m8}) into (\ref{m7}) and neglecting second order terms the system
\begin{eqnarray}
\alpha_{n+1}&=&\alpha_{n}(p+\left(\nabla
s(u^{*})\right)u^{*})+\left(\nabla s(u^{*})\right)z_{n}\nonumber\\
&=&\alpha_{n}(p+q)+\left(\nabla s(u^{*})\right)z_{n}\label{m9}\\
 z_{n+1}&=&Sz_{n},\label{m10}
\end{eqnarray}
is obtained. (Last equality in (\ref{m9}) comes from (P1) and (P2), which imply that $\left(\nabla s(u^{*})\right)u^{*}=qs(u^{*})=q$.)
Due to hypotheses (ii) and (iii), the sequence $\{z_{n}\}_{n=0}^{\infty}$ in (\ref{m10}) converges to zero. This and property (P2) imply then that $\alpha_{n}\rightarrow 0$ in (\ref{m9}), leading to local convergence. Finally, the fastest rate of convergence occurs when the factor $p+q$ in (\ref{m9}) is zero.$\Box$
}

In summary, under the hypotheses of Theorem \ref{theorem1}, the iteration map (\ref{iterop}) is contractive in a neighborhood of the fixed
point, with the fastest rate of convergence  when
$q=-p$. Condition (iii) was already obtained in
\cite{pelinovskys}, for the equations treated there and the
continuous version of the \PM method. In this sense, Theorem
\ref{theorem1} establishes the fact that (iii) is one of the sufficient
conditions for the local convergence for more general methods and in more general systems.


{
Assumption (iii) also suggests a dependence of the convergence on
the choice of the initial iterate, not only in the sense required by the local convergence, but also because $u_{0}$ must contain the correct directions. The contribution to the iteration error of the components of $u_{0}$ in these \lq harmful\rq\ eigendirections that (iii) is concerned with, can be sketched as follows. Assume for simplicity that $S$ contains one semisimple eigenvalue $\lambda_{0}$ with $|\lambda_{0}|=1$ and the rest of the spectrum (except $\lambda=p$) is below one in modulus. In (\ref{m8}), we can decompose the term $z_{n}$ in the form $z_{n}=v_{n}+w_{n}$ with $v_{n}\in Ker(\lambda_{0} I-S), w_{n}\in V\backslash Ker(\lambda_{0} I-S)$.  (Thus, $v_{n}$ can be written as a linear combination of a basis of $Ker(\lambda_{0} I-S)$, with the coordinates depending on $n$.) Now, (\ref{m9}) and (\ref{m10}) can be written as
\begin{eqnarray}
&&\alpha_{n+1}=\alpha_{n}(p+q)+\left(\nabla s(u^{*})\right)v_{n}
+\left(\nabla s(u^{*})\right)w_{n},\label{m91}\\
&& v_{n+1}=Sv_{n}=\lambda_{0}v_{n},\label{m92}\\
&& w_{n+1}=Sw_{n}.\label{m101}
\end{eqnarray}
Therefore, due to (\ref{m101}) and the previous assumptions on the spectrum of $S$, $w_{n}$ goes to zero, while (\ref{m92}) implies $v_{n}=\lambda_{0}^{n}v_{0}, n=0,1,\ldots,$ being $v_{0}$ the component of $e_{0}$ in $Ker(\lambda_{0}I-S)$ (which is to say the component of $u_{0}$ in $Ker(\lambda_{0}I-S)$). Then (\ref{m91}) becomes
 \begin{eqnarray*}
\alpha_{n+1}=\alpha_{n}(p+q)+\lambda_{0}^{n}\left(\nabla s(u^{*})\right)v_{0}
+\left(\nabla s(u^{*})\right)w_{n}.
\end{eqnarray*}
Thus, in general, $\alpha_{n}$ would be $O(||v_{0}||)$ as $n\rightarrow\infty$. As proved by Theorem \ref{theorem1}, if $v_{0}=0$ (condition (iii)) and using (P2), then $\alpha_{n}$ tends to zero as $n\rightarrow\infty$. The previous arguments also say that the errors would behave as the size of the component $v_{0}$.
}

The comparison between the matrices $S$ and $F^{\prime}(u^{*})$ reveals that the stabilizing factor acts like a filter for the harmful direction of the error that leads to the nonconvergence of the classical fixed-point algorithm in this case. The spectrum of $F^{\prime}(x^{*})$ differs from that of $S$ in the dominant eigenvalue $p$, which is transformed to some less than one (or, eventually, to zero eigenvalue if the optimal case is taken), leading to convergence if the rest of the spectrum of $S$, with probably the help of the initial iteration, behaves in the way described in Theorem \ref{theorem1} {see the numerical experiments in sections \ref{sec32} and \ref{sec34})}.


\subsection{Some examples}
\label{sec32}
As a first example, the application of the \PM method to generate
localized ground state solutions of the nonlinear Schr\"{o}dinger (NLS) model
\begin{eqnarray}
\label{doub_well11}
    iu_{t}+\partial_{xx} u+V(x)u-|u|^{2}u=0,
    \end{eqnarray}
with {a potential $V(x)$} is considered. A ground state
solution has the form $u(x,t)=e^{i\mu t}U(x)$, where $\mu\in
\mathbb{R}$ and the profile $U(x)$ is assumed to be real and
localized ($U\rightarrow 0,\; |x|\rightarrow\infty$) and then must
satisfy
\begin{eqnarray}
\label{doub_well12}
     U^{\prime\prime}(x)+V(x)U(x)-\mu U(x)-U^{3}(x)=0.
\end{eqnarray}
The \PM method (as a representative of the family (\ref{m4}), see Section \ref{sec22}) can be applied to a discretization of
(\ref{doub_well12}). One way to treat numerically the problem is
approximating (\ref{doub_well12}) on a sufficiently long interval
$(-l,l)$ and then discretizing the corresponding system for the
profile. As an illustration, the discretization based on a Fourier
collocation method for the periodic problem is taken, in such a
way that the corresponding discrete equations have the form
(\ref{m1}) with
\begin{eqnarray*}
L=D^{2}+{\rm diag}(V)-\mu I,\quad N(U_{h})=-U_{h}.^{3},
\end{eqnarray*}
where $D$ is the pseudospectral differentiation matrix,
(see \cite{Boyd}, chapter 6 and \cite{Canutohqz}, chapter 2), ${\rm diag}(V)$ is the diagonal matrix with
elements $V_{j}=V(x_{j}), x_{j}=-l+jh, j=0,\ldots,m-1$, $I$ is the
$m\times m$ identity matrix and the dot in the nonlinearity $N$
stands for the Hadamard product from the approximation $U_{h}\in
\mathbb{R}^{m}$ to the exact values of the profile at the grid
points $x_{j}$.
\begin{table}
\begin{center}
\begin{tabular}{|c|c|c|c|}
\hline\hline
\multicolumn{2}{|c|}
{$V(x)={\rm sech}^{2}(x)$} & \multicolumn{2}{|c|}{$V(x)=-6({\rm sech}^{2}(x-1)+{\rm sech}^{2}(x+1))$}\\
\multicolumn{2}{|c|}
{$\mu=1.3$} & \multicolumn{2}{|c|}{$\mu=1.43$}\\
\hline
eigs $S$&eigs $(F^{\prime}(u^{*}))$&eigs $S$&eigs $(F^{\prime}(u^{*}))$\\\hline
2.9999E+00&7.0640E-01&8.0032E+00&8.0032E+00\\
7.0640E-01&3.2731E-01&-5.6760E+00&-5.6760E+00\\
3.2731E-01&1.9060E-01&2.9999E+00&-1.5841E+00\\
1.9060E-01&1.2518E-01&-1.5841E+00&1.1350E+00\\
1.2518E-01&8.8644E-02&1.1350E+00&-9.7207E-01\\
8.8644E-02&6.6133E-02&-9.7207E-01&-5.7730E-01\\
\hline\hline
\end{tabular}
\end{center}
\caption{Ground state generation for (\ref{doub_well11})
with $V(x)={\rm sech}^{2}(x), \mu=1.3$ and $V(x)=-6({\rm sech}^{2}(x-1)+{\rm sech}^{2}(x+1)), \mu=1.43$. Six largest magnitude eigenvalues of the
approximate iteration matrix (\ref{m2b}) and of the Jacobian (\ref{m6c}). Both are evaluated
at the last computed iterate $U_{f}$ of the \PM method for (\ref{doub_well12}) and AITEM for (\ref{doub_well12}), respectively.}\label{tav1}
\end{table}

\begin{figure}[htbp]
\centering
\subfigure[]{
\includegraphics[width=6.6cm]{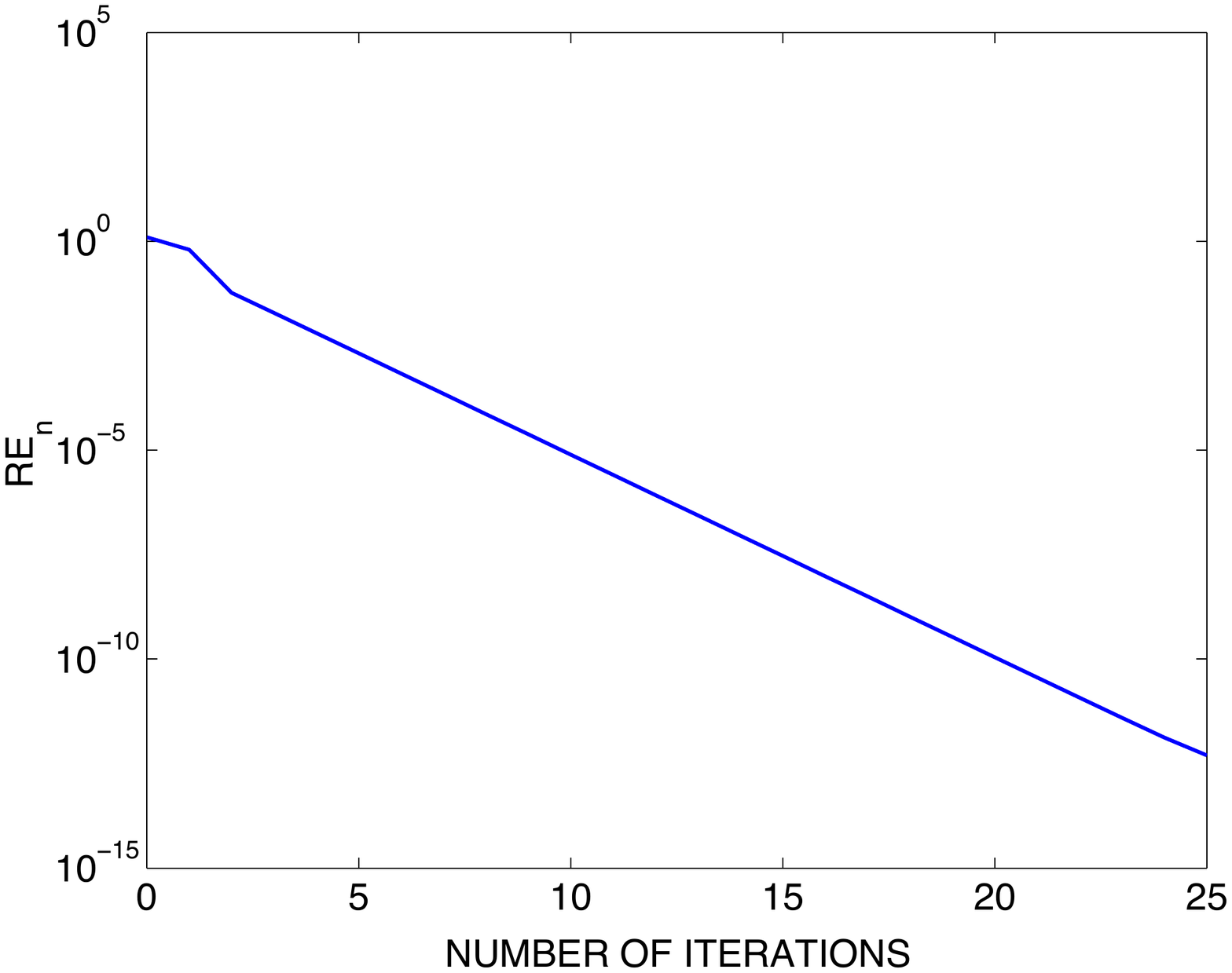} }
\subfigure[]{
\includegraphics[width=6.6cm]{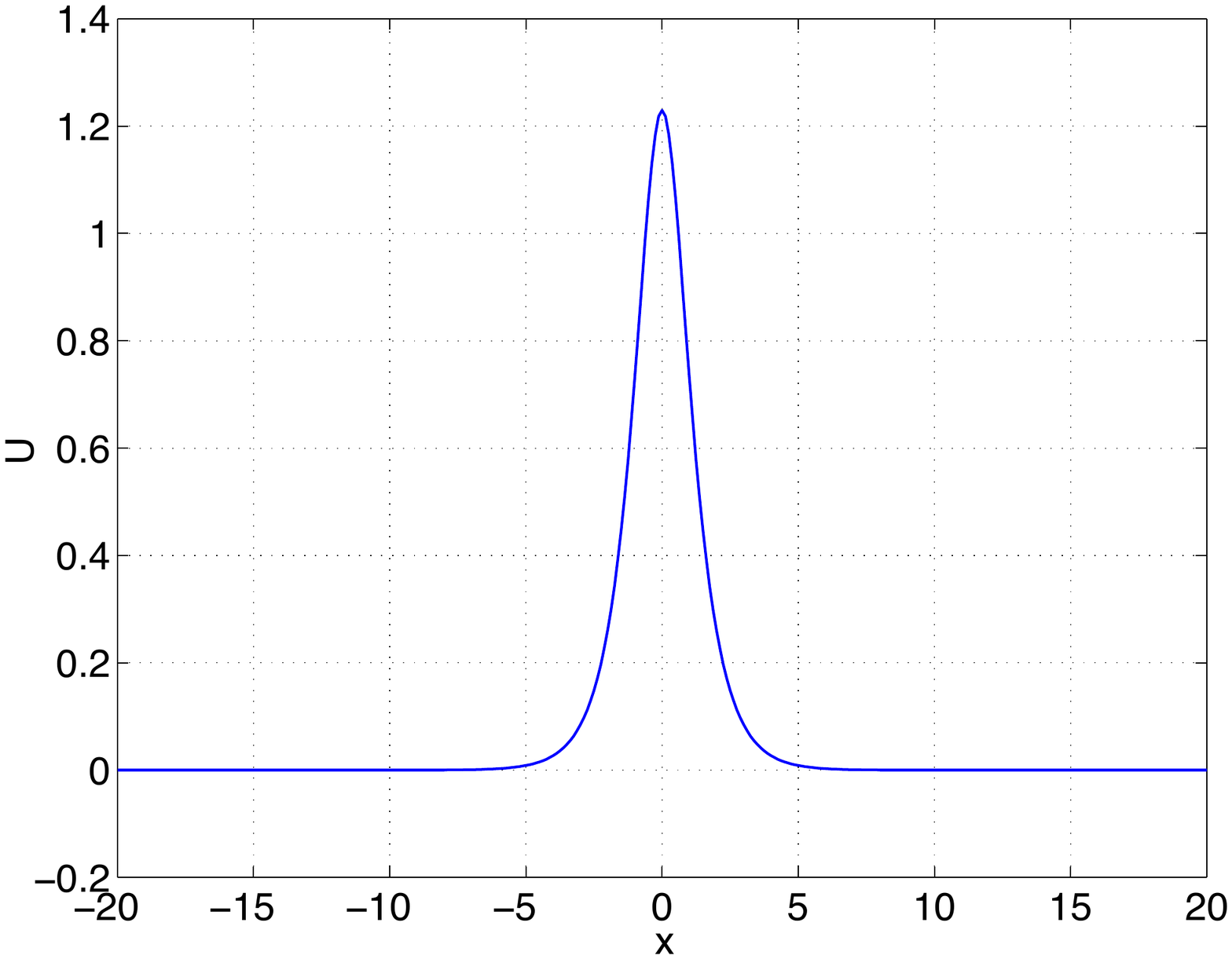} }
\caption{Convergence results of the \PM method for
(\ref{doub_well11}) and $V(x)={\rm sech}^{2}(x)$: (a)  Logarithm of
the residual errors vs number of iterations. (b) Approximate profile.} \label{fexample221}
\end{figure}

{
The ground state generation of  (\ref{doub_well11}) is illustrated with two potentials (see \cite{lakobay,yang2} and references therein for applications). The first one is $V(x)={\rm sech}^{2}(x)$.} For $\mu=1.3$ and a Gaussian profile as initial iteration, the \PM
method has been run. Table \ref{tav1} (first column) shows the six
largest magnitude eigenvalues of the approximated iteration matrix
(\ref{m2b}) at $u^{*}=U_{f}$, where $U_{f}$ is the last computed
iterate (an analytical expression for the ground state profile is
not known). The dominant eigenvalue $p=3$ (corresponding to the
degree of homogeneity for this case) is observed, with the rest
below one. The effect of the method is observed in the second column of Table \ref{tav1}, that displays the dominant eigenvalues of the Jacobian (\ref{m6c}). The magnitude of the eigenvalues is less than one, guaranteeing the convergence of the method, which is
illustrated in Figure \ref{fexample221} (a). This shows
the
logarithm of the residual error (\ref{res})
as function of the number of
iterations. In approximately $25$ iterations, a residual error of about
$1.5\times 10^{-12}$ is obtained. The ground state profile is
shown in Figure \ref{fexample221} (b). {The convergence of the stabilizing factor to one has also been checked, with a final discrepancy, in $25$ iterations,  of about $3\times 10^{-14}$.}

\begin{figure}[htbp]
\centering
\includegraphics[width=0.7\textwidth]{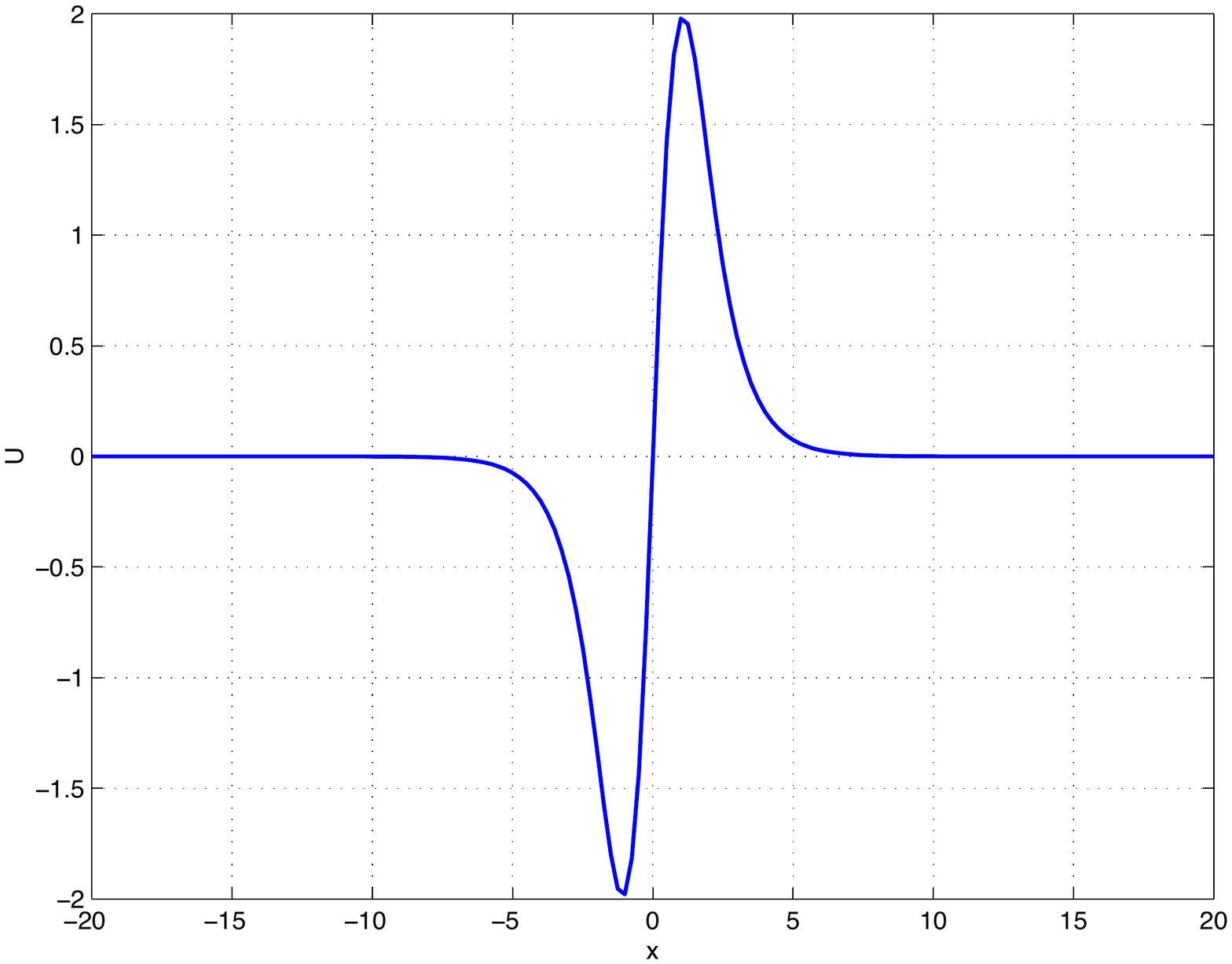}
\caption{Antisymmetric solution of (\ref{doub_well12}): numerical
profile obtained with AITEM, $\mu=1.43$.} \label{fexample222}
\end{figure}

{
As a second example, we consider (\ref{doub_well11}) with a double-well
potential
 $V(x)=-6({\rm sech}^{2}(x-1)+{\rm sech}^{2}(x+1))$.
As indicated in \cite{lakobay}, the \PM method fails in the search
for an anti-symmetric solution of (\ref{doub_well12}) (see the profile in Figure \ref{fexample222}, obtained for $\mu=1.43$ with the AITEM method, \cite{yangl1}). This case of divergence can be justified using the previous results.}
For $\mu=1.43$,
Table \ref{tav1} (third column) shows the six largest magnitude
eigenvalues of the iteration matrix evaluated at the profile
obtained with the AITEM method. Besides the
eigenvalue $\lambda=3$, associated to the degree of homogeneity of
the nonlinear part, Table \ref{tav1} reveals the existence of
other eigenvalues with magnitude above one. This divergence is also confirmed by the eigenvalues of the Jacobian (\ref{m6c}), shown in the fourth column of Table \ref{tav1}.

\subsection{Systems with symmetries. Orbital convergence}
\label{sec33}
In many situations, the system (\ref{m1}) admits symmetries, \cite{olver,marsdenr}. This means that there is an ${\nu}$-parameter group of transformations ($\nu\geq 1$)
\begin{eqnarray}
\mathcal{G}=\{G_{\alpha}:\mathbb{R}^{m}\rightarrow\mathbb{R}^{m},
\alpha=(\alpha_{1},\ldots,\alpha_{\nu})\in\mathbb{R}^{\nu}\},\label{m11}
\end{eqnarray}
with the property of transforming solutions of (\ref{m1}) into other solutions:
\begin{eqnarray}
Lu^{*}=N(u^{*})\Rightarrow L(G_{\alpha}u^{*})=N(G_{\alpha}u^{*}),\quad
\alpha\in\mathbb{R}^{\nu}.\label{m12}
\end{eqnarray}
For simplicity, we assume that the transformations $G_{\alpha}$ in (\ref{m11}) are smooth and $\mathcal{G}$ is Abelian. {(In traveling wave generation, typical examples are, as mentioned before, translations and phase rotations, see Section \ref{sec34}.)}. The existence of a symmetry group for (\ref{m1}) has several consequences. We emphasize two of them:
\begin{itemize}
\item The group (\ref{m11}) defines orbits of solutions of (\ref{m1}):
\begin{eqnarray*}
\mathcal{G}(u^{*})=\{G_{\alpha}u^{*}:\alpha\in\mathbb{R}^{\nu}\}.
\end{eqnarray*}
The space of solutions of (\ref{m1}) is partitioned into these orbits, in such a way that
assumption (H1) in Theorem \ref{theorem1} does not hold in this case: a fixed point $u^{*}$ cannot be isolated and the concept of convergence for the iterative methods must be redefined. Under these conditions, it is said that the iteration (\ref{m4}) is orbitally convergent to $u^{*}$ if $u_{n}$ converges to $G_{\alpha}u^{*}$ for some $\alpha\in \mathbb{R}^{\nu}$.
\item The pencil $A(\lambda)$ admits $\lambda=1$ as eigenvalue, since differentiation with respect to $\alpha$ in (\ref{m12}) implies
\begin{eqnarray*}
\left(L-N^{\prime}(u^{*})\right)\frac{\partial}{\partial
\alpha_{j}}\big|_{\alpha=0}G_{\alpha}(u^{*})=0,\quad
j=1,\ldots, \nu,
\end{eqnarray*}
and the infinitesimal generators of the group, \cite{olver}
\begin{eqnarray}
u\mapsto v_{j}(u)=\frac{\partial}{\partial
\alpha_{j}}\big|_{\alpha=0}G_{\alpha}(u),\quad
j=1,\ldots, \nu,\label{igg}
\end{eqnarray}
evaluated at $u=u^{*}$, are associated eigenvectors.
\end{itemize}
{The convergence result in Theorem \ref{theorem1} can be adapted to this case as follows. First, hypothesis (H1) is substituted by
\begin{itemize}
\item[(H1)'] $\mathcal{G}$ is a symmetry group of (\ref{m1}) with $dim Ker A(1)=\nu$.
\end{itemize}
As far as the spectrum of $A(\lambda)$ is concerned, we still assume (i) and (ii) of Theorem \ref{theorem1} while the third condition is now
\begin{itemize}
\item[(iii)] If $|\lambda|=1\Rightarrow\left\{\begin{matrix}\lambda \quad \mbox {is
semisimple}\\\mbox{if} \quad \lambda\neq 1\Rightarrow
u_{0} \quad \mbox {does not have  } \\
\mbox {component in}\quad Ker A(\lambda)\end{matrix}\right.$
\end{itemize}
In this case the convergence is orbital, in the sense above defined. Note that now the errors $e_{n}, n=0,1,\ldots,$ can be decomposed in the form (cf. (\ref{m8}))
\begin{eqnarray*}
e_{n}=\alpha_{n}u^{*}+\sum_{k=1}^{\nu}\beta_{n,k}v_{k}(u^{*})+z_{n}, \quad
\alpha_{n},\beta_{n,k}\in \mathbb{R}, k=1,\ldots,\nu,\quad z_{n}\in V,
\end{eqnarray*}
where $\{v_{k}(u^{*})\}_{k=1}^{\nu}$ is the basis (\ref{igg}) of $Ker A(1)$ and now $V$ is the (unique) supplementary $(m-\nu-1)$ dimensional space of $span(u^{*})+Ker A(1)$
with $S(V)\subset V$. Now, the
sequences $\alpha_{n}, \beta_{n,k}, =1,\ldots,\nu$ and $z_{n}$ satisfy
\begin{eqnarray}
\alpha_{n+1}&=&\alpha_{n}(p+q)+\sum_{k=1}^{\nu}\beta_{n,k}\left(\nabla s(u^{\ast})\right)v_{k}(u^{*})+\left(\nabla s(u^{\ast})\right)z_{n},\label{mb1}\\
z_{n+1}&=&Sz_{n},\label{mb2}\\
\beta_{n+1,k}&=&\beta_{n,k}, \quad k=1,\ldots,\nu.\label{mb3}
\end{eqnarray}
Note, on the other hand, that (P1) and (\ref{m12}) imply
\begin{eqnarray}
s(G_{\alpha}(u^{*}))=1,\quad \alpha=(\alpha_{1},\ldots,\alpha_{\nu})\in\mathbb{R}^{\nu}.\label{ese}
\end{eqnarray}
Differentiating (\ref{ese}) with respect to each $\alpha_{j}, j=1,\ldots,\nu$ and evaluating at $\alpha=0$ we have
\begin{eqnarray*}
\left(\nabla s(u^{\ast})\right)v_{j}=0,\quad j=1,\ldots,\nu,
\end{eqnarray*}
and (\ref{mb1})-(\ref{mb3}) can be written as
\begin{eqnarray*}
\alpha_{n+1}&=&\alpha_{n}(p+q)+\left(\nabla s(u^{\ast})\right)z_{n},\\
z_{n+1}&=&Sz_{n},\\
\beta_{n+1,k}&=&\beta_{n,k}, \quad k=1,\ldots,\nu\\
&&\Rightarrow \beta_{n,k}=\beta_{0,k},\quad  k=1,\ldots,\nu, \quad n\geq 0.
\end{eqnarray*}
Consequently,
\begin{eqnarray*}
G_{(\beta_{0,1},\ldots,\beta_{0,\nu})}(u^{*})+\alpha_{n}x^{*}+z_{n}
\end{eqnarray*}
differs from
\begin{eqnarray*}
u_{n}=u^{*}+e_{n}=u^{*}+\sum_{k=1}^{\nu}\beta_{0,k}v_{k}(u^{*})+z_{n},
\end{eqnarray*}
in $O(||e_{n}||^{2})$ terms. Under the above mentioned hypotheses,
the convergence is to the
element $G_{(\beta_{0,1},\ldots,\beta_{0,\nu})}(u^{*})$ of the orbit of $u^{*}$,
determined by the component of the initial iteration in
Ker$(I-S)$.As in Theorem \ref{theorem1}, the fastest rate of convergence occurs when $q=-p$.
}

\subsection{Some examples}
\label{sec34}
This case is illustrated by the generation of solitary wave solutions of nonlinear Schr\"{o}dinger equations of the form (see e.~g. \cite{sulems} and references therein)
\begin{eqnarray}
iu_{t}+u_{xx}+|u|^{2\sigma}u=0,\quad -\infty<x<\infty,\quad t>0,\label{nlse1}
\end{eqnarray}
where $\sigma>0$. The symmetry group for (\ref{nlse1}) consists of gauge transformations and translations
\begin{eqnarray}
G_{{\theta_{0},x_{0}}}(u(x))=e^{i\theta_{0}}u(x+x_{0}),\quad \theta_{0},x_{0}\in \mathbb{R}.\label{nlse2}
\end{eqnarray}
Solitary wave solutions of (\ref{nlse1}) can be obtained from profiles {$U(x)$ satisfying}
\begin{eqnarray}
U^{\prime\prime}+|U|^{2\sigma}u-\lambda_{1}U-i\lambda_{2}U^{\prime}=0,\label{nlse3}
\end{eqnarray}
for some real parameters $\lambda_{1},\lambda_{2}$. This leads to the explicit formulas
\begin{eqnarray}
U(x)&=&\rho(x)e^{i\theta(x)}\label{nlse4}\\
\rho(x)&=&(a(\sigma+1))^{1/2\sigma}\left({\rm sech}(\sigma\sqrt{a}x)\right)^{1/\sigma},\quad a=\lambda_{1}-(\lambda_{2}^2)/4,\label{nlse4b}\\
\theta(x)&=&\frac{\lambda_{2}}{2}x.\label{nlse4c}
\end{eqnarray}
Due to the symmetry group (\ref{nlse2}), the two-parameter orbit of the solution given by (\ref{nlse4})-(\ref{nlse4c}) is of the form
\begin{eqnarray}
\mathcal{G}(\rho,\theta)=\{\varphi=\rho(x-x_{0})e^{i\theta(x-x_{0})+i\theta_{0}}: x_{0},\theta_{0}\in \mathbb{R}\}.\label{nlse41}
\end{eqnarray}
The four-parameter family of solitary wave solutions of (\ref{nlse1}) is finally of the form
\begin{eqnarray*}
\psi(x,t,a,c,x_{0},\theta_{0})=G_{(t\lambda_{1},t\lambda_{2})}(\varphi)=\rho(x-ctx_{0})e^{i\theta(x-ct-x_{0})+i\theta_{0}+i(a+(c^{2}/4))t},\label{nlse5}
\end{eqnarray*}
(where $c=\lambda_{1}$). As far as the discretization is concerned, the corresponding Fourier collocation approximation of (\ref{nlse3})
\begin{eqnarray*}
D^{2}U_{h}+|(U_{h})|.^{2\sigma}.U_{h}-\lambda_{1}U_{h}-i\lambda_{2}DU_{h}=0,
\end{eqnarray*}
(the dot stands for the Hadamard product)
inherites the symmetry group infinitesimally generated by (see (\ref{igg}))
\begin{eqnarray*}
U_{h}\mapsto v_{1}(U_{h})=iU_{h},\quad U_{h}\mapsto v_{2}(U_{h})=DU_{h}.
\end{eqnarray*}
They are associated, respectively, to phase rotations and spatial translations.

\begin{table}
\begin{center}
\begin{tabular}{|c|c|}
\hline\hline
$\sigma=1$ & $\sigma=2$\\
\hline
2.9999E+00&4.9999E+00\\
9.9999E-01&9.9999E-01\\
9.9999E-01&9.9999E-01\\
4.9999E-01&4.2857E-01\\
3.3333E-01&2.3810E-01\\
2.9999E-01&1.9999E-01\\
\hline\hline
\end{tabular}
\end{center}
\caption{Solitary wave generation of  (\ref{nlse1}). Six largest magnitude eigenvalues of the
approximate iteration matrix (\ref{m2b}) evaluated
at the exact profile (\ref{nlse4}) with $\lambda_{1}=\lambda_{2}=1, x_{0}=\theta_{0}=0$ for $\sigma=1$ (first column) and $\sigma=2$ (second column).}\label{tav2}
\end{table}

Local convergence of the \PM method is first checked by Table \ref{tav2}. This shows, for $\sigma=1,2$, the first largest magnitude eigenvalues of the iteration matrix $S$ at the exact profile (\ref{nlse4})-(\ref{nlse4c}) with parameters $\lambda_{1}=\lambda_{2}=1$, (where $x_{0}=\theta_{0}=0$). The results guarantee the satisfaction of the conditions for local convergence. (Note that, in this case, since the symmetry group is two dimensional, the eigenvalue $\lambda=1$ has geometric multiplicity equals two.)

\begin{figure}[htbp]
\centering \subfigure[]{
\includegraphics[width=6.6cm]{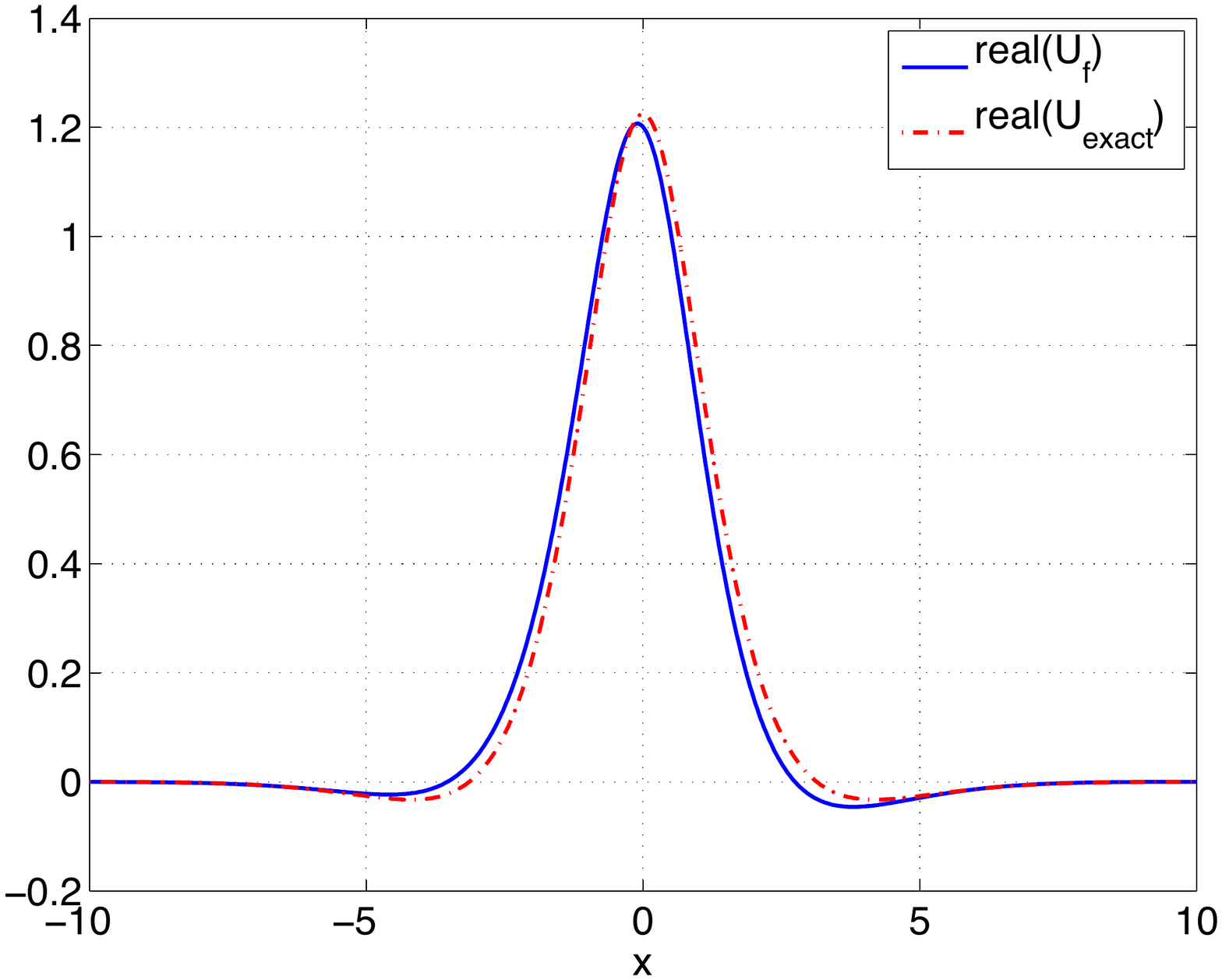} }
\subfigure[]{
\includegraphics[width=6.6cm]{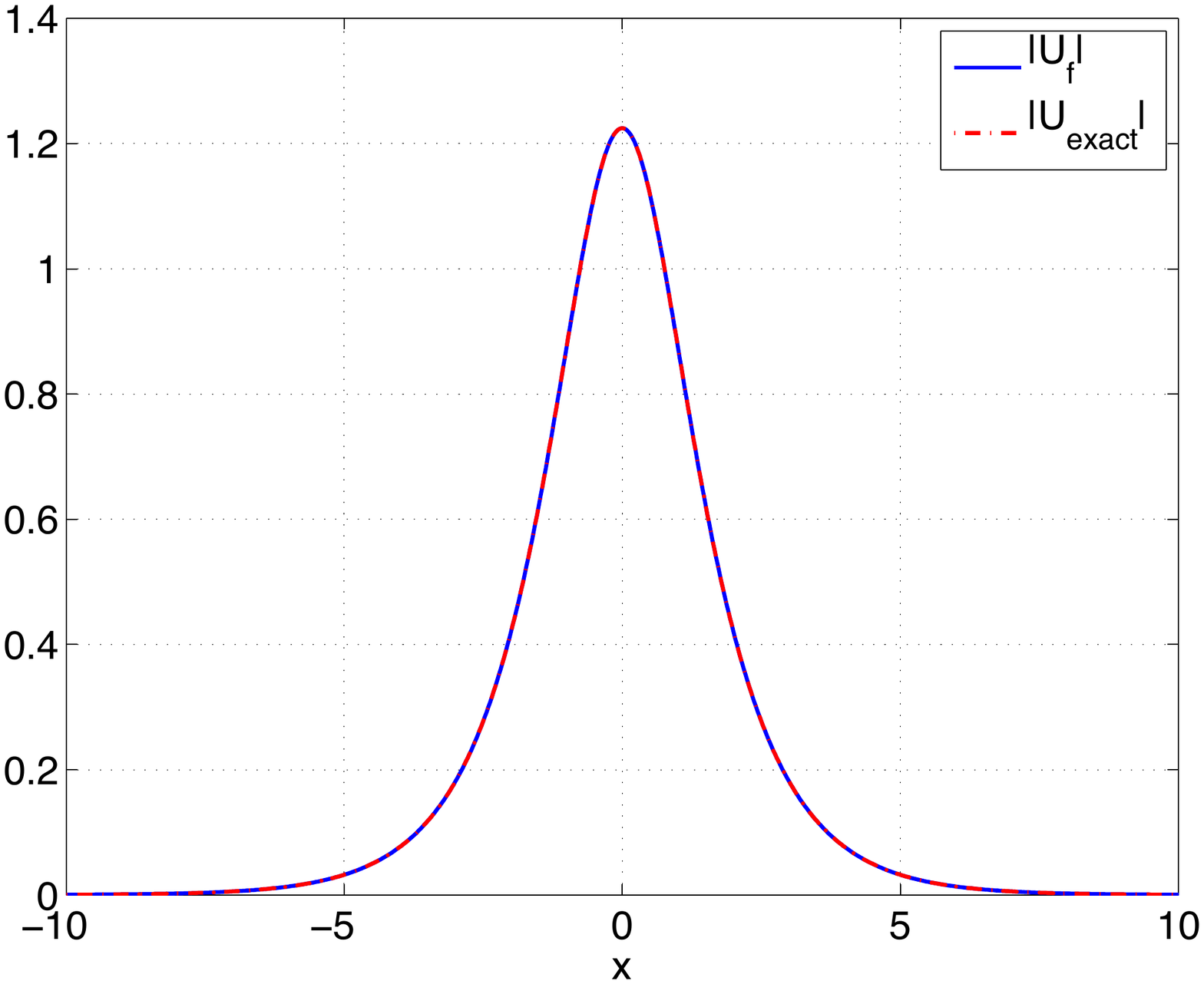} }
\subfigure[]{
\includegraphics[width=6.6cm]{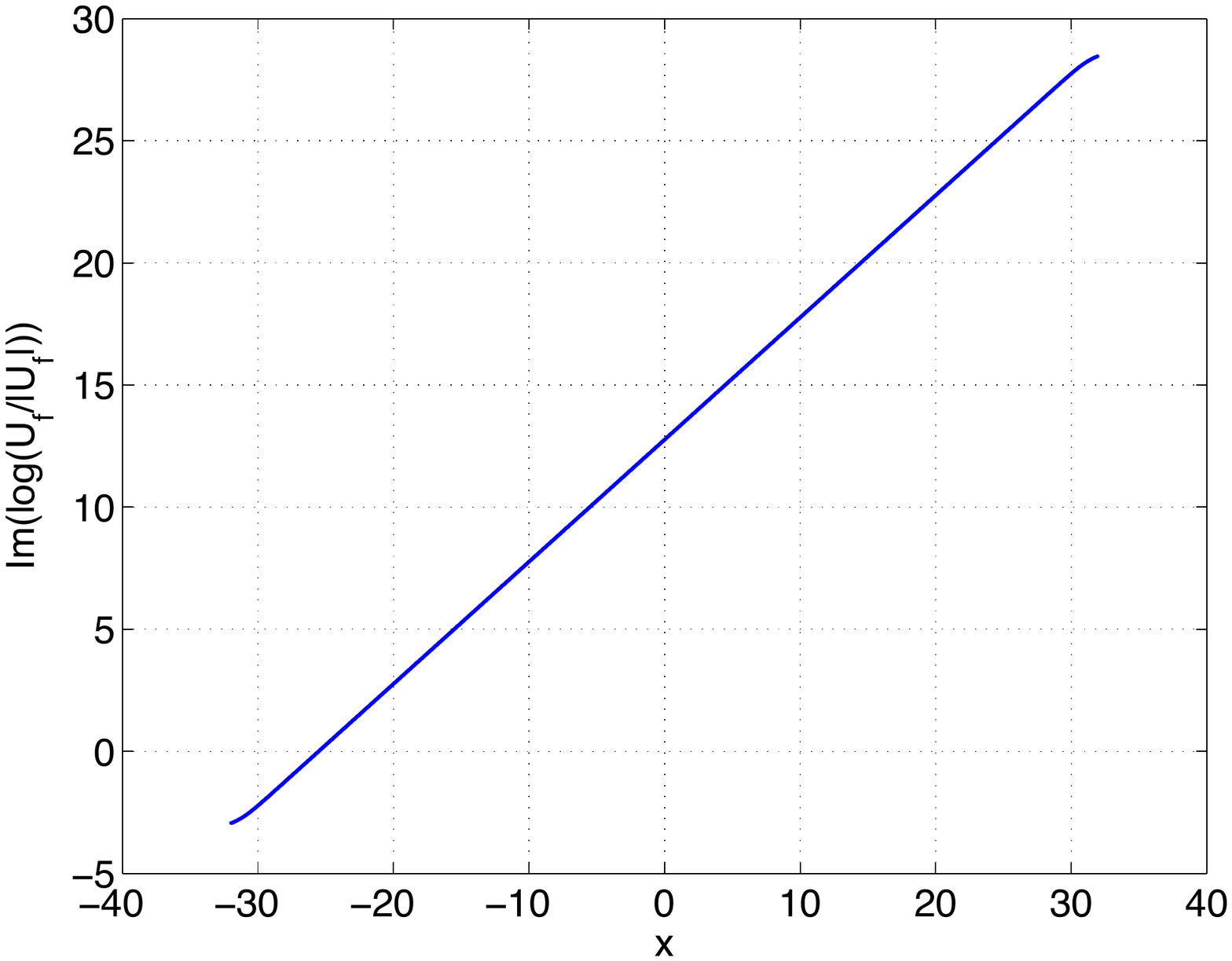} }
\caption{{Solitary wave generation of (\ref{nlse1}). The parameters are $\lambda_{1}=\lambda_{2}=1, x_{0}=\theta_{0}=0, \epsilon_{1}=0.2, \epsilon_{2}=0$. (a) Real part of the last computed iterate $U_{f}$ obtained by the Petviashvili scheme (solid line) and of the exact profile $U_{\rm exact}$ (dashed line). (b) Modulus of the last computed iterate $U_{f}$ obtained by the Petviashvili scheme (solid line) and of the exact profile $U_{\rm exact}$ (dashed line). (c) Fitting line to the phase of the computed profile $U_{f}$.} \label{fnls1}}
\end{figure}

{
The following results illustrate the orbital convergence. An exact profile of the form (\ref{nlse4})-(\ref{nlse4c}) with $\lambda_{1}=\lambda_{2}=1, x_{0}=\theta_{0}=0$, denoted by $U_{\rm exact}$, is perturbed in the form
\begin{eqnarray}
U_{0}=U_{\rm exact}+\epsilon_{1} i U_{\rm exact}+\epsilon_{2} D U_{\rm exact},\label{m36}
\end{eqnarray}
(with small parameters $\epsilon_{1},\epsilon_{2}$). Now the \PM method is run with (\ref{m36}) as initial iteration. The results are illustrated by two experiments. Figure \ref{fnls1}(a) (resp. Figure \ref{fnls1}(b)) compares the real part (resp. the modulus) of the exact profile $U_{\rm exact}$ with that of the last computed iterate, denoted by $U_{f}$ and obtained with a residual error below $10^{-13}$ and for $\epsilon_{1}=0.2, \epsilon_{2}=0$.
}

{While moduli are practically indistinguishable, Figure \ref{fnls1}(a) reveals a phase displacement of the computed profile with respect to the exact one. The phase of $U_{f}$ (computed as $Im\left(log\left(U_{f}/|U_{f}|\right)\right)$, modulo $2\pi$) has been calculated. The resulting data are fitted to a line $y=mx+n$, see Figure \ref{fnls1}(c). The computed slope is $m=4.9934\times 10^{-1}$ (approximating the corresponding value $\lambda_{2}/2$ of (\ref{nlse4})-(\ref{nlse4c}) for this case) while $n=1.2764\times 10^{1}$, which modulo $2\pi$ is $1.9751\times 10^{-1}$, an approximation to the value of $\epsilon_{1}$. These results suggest that the computed profile is closer to the element of the orbit (\ref{nlse41}) of $U_{\rm exact}$ with new phase $\theta_{0}+\epsilon_{1}=\epsilon_{1}$ and the same translational parameter $x_{0}=0$.
}

\begin{figure}[htbp]
\centering \subfigure[]{
\includegraphics[width=6.6cm]{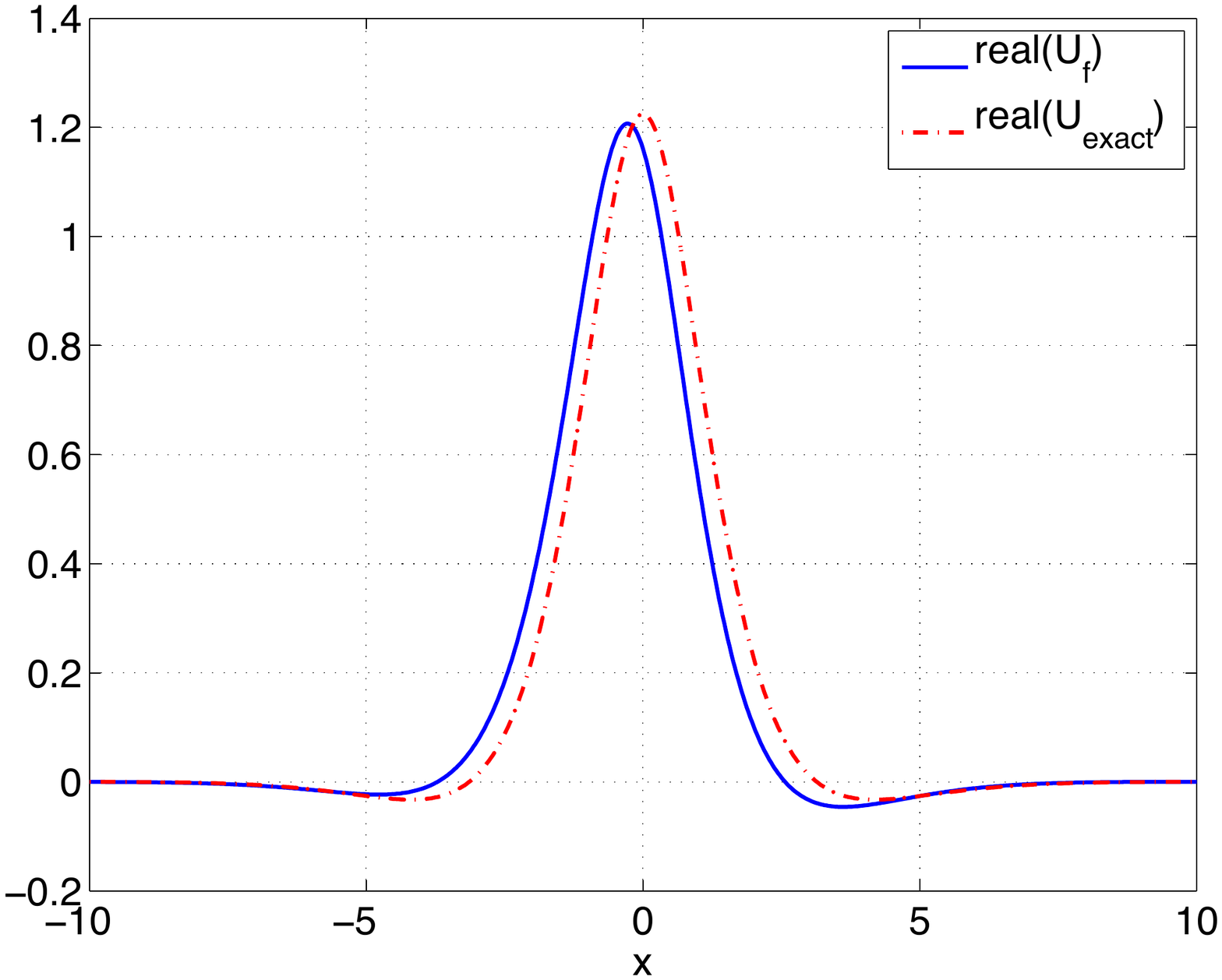} }
\subfigure[]{
\includegraphics[width=6.6cm]{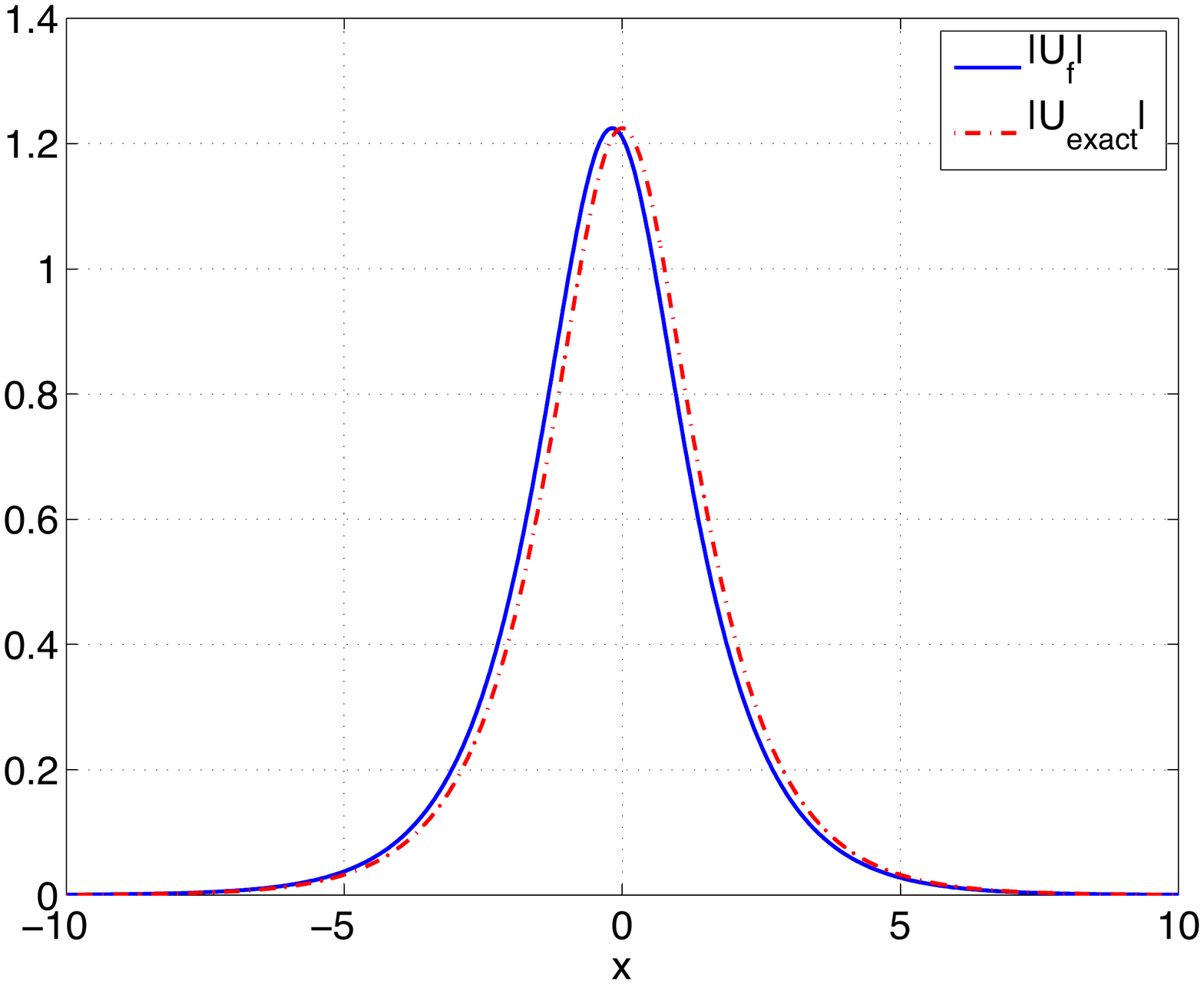} }
\subfigure[]{
\includegraphics[width=6.6cm]{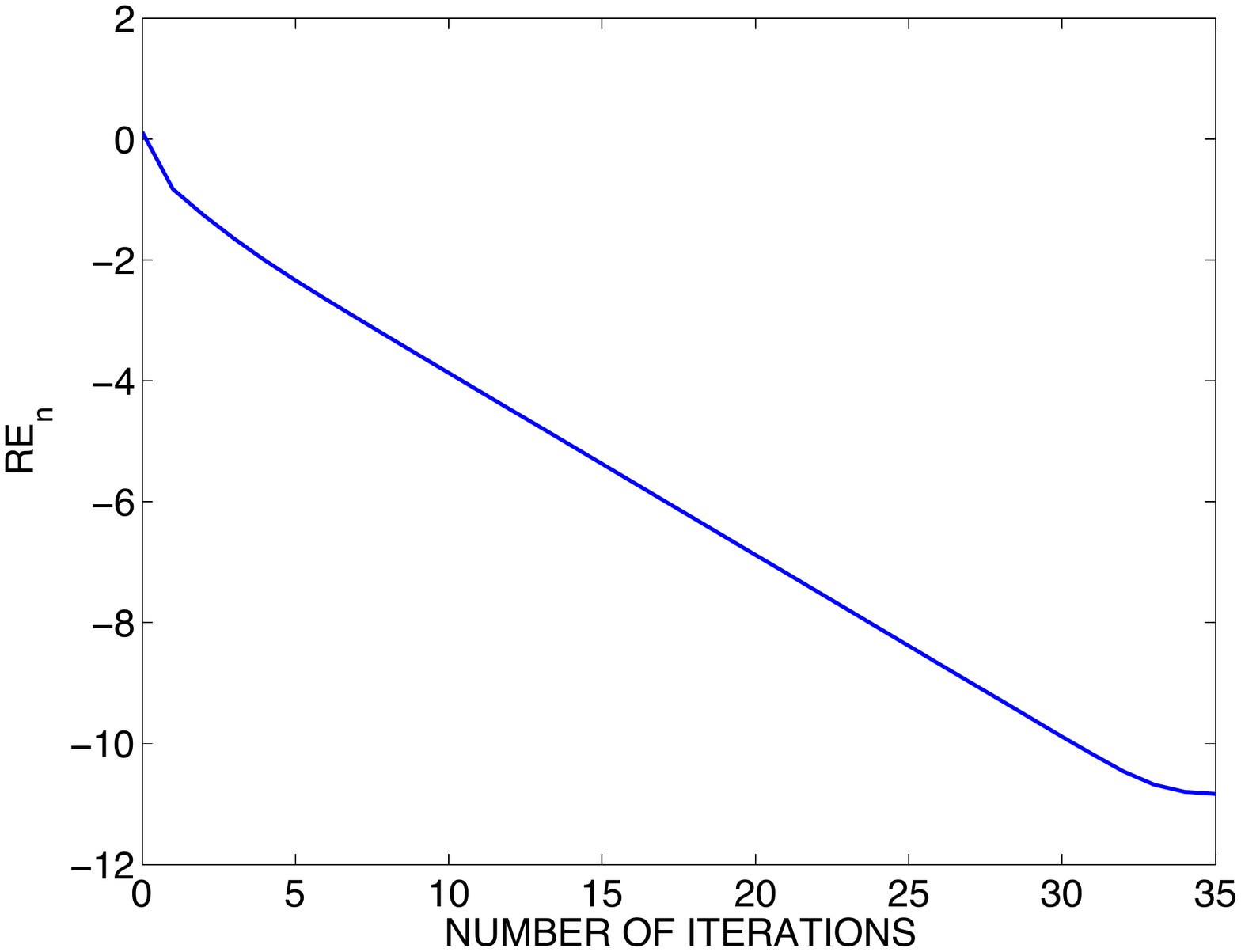} }
\caption{{Solitary wave generation of (\ref{nlse1}). The parameters are $\lambda_{1}=\lambda_{2}=1, x_{0}=\theta_{0}=0, \epsilon_{1}=0.2, \epsilon_{2}=0.2$. (a) Real part of the last computed iterate $U_{f}$ obtained by the Petviashvili scheme (solid line) and of the exact profile $U_{\rm exact}$ (dashed line). (b) Modulus of the last computed iterate $U_{f}$ obtained by the Petviashvili scheme (solid line) and of the exact profile $U_{\rm exact}$ (dashed line). (c)  Logarithm of the residual error against number of iterations.}} \label{fnls2}
\end{figure}

{A second experiment is performed with the values $\epsilon_{1}=0.2, \epsilon_{2}=0.2$. Figure \ref{fnls2}(a) (resp. Figure \ref{fnls2}(b)) compares the real part (resp. the modulus) of the exact profile $U_{\rm exact}$ with that of the last computed iterate $U_{f}$, after $35$ iterations and with a residual error (\ref{res}) below $10^{-11}$ (see Figure \ref{fnls2}(c)). The error in the stabilizing factor is of the order of the machine precision (below $10^{-15}$). In this case, besides the phase shift, also the modulus is affected by a displacement (cf. Figure \ref{fnls1}(b)). The corresponding fitting line to the phase data of the computed profile $U_{f}$, $y=mx+n$ has a slope $m=4.9933\times 10^{-1}$ while $n=1.2852\times 10^{1}$, which modulo $2\pi$ is approximately $2.8553\times 10^{-1}$. This last value approximates $\epsilon_{1}+(\lambda_{2}/2)\epsilon_{2}=0.3$, suggesting that the computed profile is close to the exact one of the form (\ref{nlse4}), (\ref{nlse41}) with group parameters $\theta_{0}+\epsilon_{1}, x_{0}+\epsilon_{2}$.
}

\section*{Acknowledgements}
This research has been supported by  projects
MTM2010-19510/MTM (MCIN), MTM2012-30860(MECC) and VA118A12-1 (JCYL). The authors want to thank the reviewers for their fruitful suggestions and comments.

\end{document}